 \documentclass[12pt,leqno]{amsart}
\newtheorem{dummy}{}[section]
\newtheorem{definition}[dummy]{Definition}
\newtheorem{theorem}[dummy]{Theorem}

\newtheorem{proposition}[dummy]{Proposition}
\newtheorem{lemma}[dummy]{Lemma}

\newtheorem{remark}[dummy]{Remark}
\newtheorem{example}[dummy]{Example}
\newtheorem{question}[dummy]{Question}
\begin{document}
\bibliographystyle{plain}
\title{Derived Langlands II: Hyper-Hecke algebras, monocentric relations and
${\mathcal M}_{cmc,\underline{\phi}}(G)$-admissibility }
\author{Victor P. Snaith}
\date{9 September 2019}
\maketitle
 \tableofcontents  
 
 \begin{abstract}
 In the book ``Derived Langlands'' I claimed to give a construction which enabled one to consider admissible representations of a locally profinite Lie group as chain complexes of monomial representations - i.e. representations given by compactly supported induction from one-dimensional representations of 
 compact, open modulo the centre subgroups. In the jargon, this process amounts to embedding the smooth representation category into a derived category of the monomial category. The book went on to illustrate some aspects of the construction as they pertained to different parts of the Langlands Programme. 

The generality described above may appear rather daunting. However, the constructions are potentially interesting and useful even in the simple case of finite-dimensional representations of finite groups or Galois groups - particularly over fields of positive characteristic. Therefore I shall spend the first two lectures introducing the monomial category in this case, where it was developed by Robert Boltje. The morphisms form a ring, which I call the hyperHecke algebra. I shall introduce the notion of exactness in this category which is necessary to define  monomial resolution of a representation. I shall spend some time on the case of complex representations of general linear groups over finite fields which illustrates several features - invariants of representations such as Kondo-Gauss sums, bijective correspondences such as Shintani base change - which are simple analogues of important ingredients of the Langlands Programme.

Other sections will describe (i) how the monomial resolution construction and the hyperHecke algebra  generalises to the case of admissible representations of compact, open modulo the centre subgroups
and (ii) how algebraic topology extends the  resolutions to admissible representations of an arbitrary locally $p$-adic Lie group. The circumstances under which this works is that of ``monomial admissibility''  
which bears the same relation to the hyperHecke algebra as does the (weaker) usual notion of admissibility to the usual Hecke algebra. 

There are several other bits and pieces: (i) \ how the famous Bernstein centre can be approached via the hyperHecke algebra (Deligne's essay on the subject uses the Hecke algebra but rather differently)
(ii) \  how the morphisms in the monomial category may be described by means of the convolution product in the classical Hecke algebra (iii) \ monomial resolutions of automorphic representations and the way spaces of classical modular forms appear in the notion of monomial exactness.

\end{abstract}

 \section{Introduction}
 
  In the book ``Derived Langlands'' (\cite{Sn18}) I claimed to give a construction which enabled one to consider admissible representations of a locally profinite Lie group as chain complexes of monomial representations - i.e. representations given by compactly supported induction from one-dimensional representations of 
 compact, open modulo the centre subgroups. In the jargon, this process amounts to embedding the smooth representation category into a derived category of the monomial category. The book went on to illustrate some aspects of the construction as they pertained to different parts of the Langlands Programme. 

The generality described above may appear rather daunting. However, the constructions are potentially interesting and useful even in the simple case of finite-dimensional representations of finite groups or Galois groups - particularly over fields of positive characteristic. 
 
  This essay is a very intuitive development of the ideas in \cite{Sn18}. The book was little more than a curiosity for the experts in the Langlands Programme\footnote{For example the Lecture Notes Series of the L.M.Soc. took several years to realise that no referee could be found for it! Expert email comments to me were not very encouraging (from P.S.``Why are you doing this?'' (curiosity!) ; from G.H. a cheerful joke about New Year resolutions (Happy New Year!)). Eventually I was approached by World Scientific, who were so helpful that I was delighted to take up their offer of publication.}. In this context ``intuitive'' is intended to indicate that this essay has the sort of tangential rigour-level which is only to be expected in the author's context\footnote{ Namely: age, retirement, mathematical isolation which is so customary in the UK, health \cite{JMOS18}, dilettante ignorance! In addition, having percolated these ideas since 1985 - and pursued them with more apparent progress since 2005 - circumstances have encouraged me to develop this project at a leisurely pace.}.
  
A provisional version of this essay was up-loaded to my home-page last year in order to be there in time for a lecture series given at the University of Sheffield in September 2019. I shall now outline the contents together with related remarks about my intentions with regard to further developments which are envisioned for``Derived Langlands III''.

Section Two introduces the hyperHecke algebra, which is the algebra of morphisms in the monomial category associated to a locally profinite group $G$ and a fixed central character $\underline{\phi}$. The hyperHecke algebra originated in equivariant stable homotopy associated to the Segal conjecture (e.g. \cite{MSZ89}). The stable homotopy setting was enough to derive the original formula for explicit Brauer  induction, taking place in the Grothendieck group of objects in the monomial category, after a suitable algebraic completion \cite{Sn89}. The completion was removed in \cite{Sn88} and the topic was elaborated upon in \cite{Sn94}. The hyperHecke algebra is given by generators and relations which are realised as morphisms between compactly induced representations obtained from a subgroup $H$ containing the centre of $G$, compact modulo the centre, and a character $\phi$ of $H$ extending the central character. When $\phi$ is trivial one gets essentially the Hecke algebra - hence the ``hyper''.
The classical Hecke algebra of $G$ is an idempotented convolution algebra. The hyperHecke algebra is also idempotented (see Section Six) and in Section Seven I give the formula for monomial morphisms in terms of the convolution algebra - more or less completely when the $\phi$'s applied to a compact open $H$ have finite image\footnote{Characteristically I have left a few details as cliff-hangers, with a remark or two intended to finish off the formula in general eventually. Ho-hum...!!}.

 Section Three calculates the centre of the monomial category.  This is motivated by attempting to approach the famous Bernstein centre of the category of admissible representations of a locally profinite group $G$ (\cite{BZ76}, \cite{BZ77}). The connection goes as follows: the monocentre gives canonical self-maps of the monomial resolutions of a representation and therefore a morphism in the Bernstein centre, as explained in Section Eleven. An element of the monocentre gives a morphism of the whole resolution - which is too strong a condition - since the morphism on the representation is induced by the first few morphisms in the  
 resolution. However, it is not unreasonable to study the Bernstein centre via the hyperHecke algebra in view of the observation that Deligne's essay on the subject \cite{PD84} uses the Hecke algebra (but rather differently). Some of the background for \cite{PD84}  is collected in Appendix Section Twelve.
 
Section Three fosters a remark about Langlands quotient theorem \cite{RPL} and Silberger's extension to the case of $p$-adic reductive groups \cite{AJS78}. These results classify irreducible admissible representation by proving that they occur as the unique irreducible quotients of compactly supported induced admissible irreducible representations suitably chosen $p$-pairs. From the monomial resolution point of view, if we are in the situation in which Schur's Lemma holds, which is usually the case, that analysing whether a morphism (or partial morphism) of the resolution induces a scalar endomorphism on the resolved representation turns the quotient theorem into a rather fancy extension of the Double Coset Formula\footnote{Progress with this must wait for ``Derived Langlands'' III}.
 
 Section Four introduces the rather obvious extension of smoothness and admissibility to which one is led once the characters $\phi$ are taken seriously - as in the hyperHecke algebra. In the classical case I show that this notion of admissibility coincides classical one. I include a few questions and remarks about the relevance of this notion to the di-p-adic situation\footnote{This is called ``the $p$-adic Langlands programme'' by the  experts, I believe.}. On this type of admissibility, I owe the reader the details analogous to the classical notes by Bill Casselman \cite{WCnotes} and Jacquet's theory \cite{AJS79}, which is intended to be given in ``Derived Langlands'' III.
 
 Sections Eight, Nine and Ten construct the bar-monomial resolution  in a manner that differs slightly from my original in \cite{Sn18}. The difference is to emphasis that, by means of a functor from the monomial category to the category of modules over the hyperHecke algebra, the bar-monomial resolution may be identified with the bar resolution over the latter.
 
 Appendix Section Twelve contains the dictionary for passing from the function space model for an induced representation to the tensor product model.

 \section{The hyperHecke algebra }
 
 Let $G$ be a locally profinite group and let $k$ be an algebraically closed field. Suppose that $\underline{\phi} : Z(G) \longrightarrow k^{*}$ is a fixed $k$-valued, continuous character on the centre $Z(G)$ of $G$. Let ${\mathcal M}_{cmc, \underline{\phi}}(G)$ be the poset of pairs $(H, \phi)$ where $H$ is a subgroup of $G$, containing $Z(G)$, which is compact, open modulo the centre of $G$ and $\phi : H \longrightarrow k^{*}$ is a $k$-valued, continuous character whose restriction to $Z(G)$ is $\underline{\phi}$.
 
 We define the hyperHecke algebra, ${\mathcal H}_{cmc}(G)$,  to be the $k$-algebra given by the following generators and relations. For $(H, \phi), (K, \psi) \in {\mathcal M}_{cmc,\underline{\phi}}(G)$,  write 
$[(K, \psi), g, (H, \phi)]$ for any triple consisting of $g \in G$, characters $\phi, \psi$ on 
subgroups $H, K \leq G$, respectively such that 
\[    (K, \psi) \leq (g^{-1}Hg, (g)^{*}(\phi)) \]
which means that $K \leq  g^{-1}Hg$ and that $\psi(k) = \phi(h)$ where $k = g^{-1}hg$ for
 $h \in H, k \in K$.

 Let ${\mathcal H}$ denote the $k$-vector space with basis given by these triples. Define a product on these triples by the formula 
\[  [(H, \phi), g_{1}, (J, \mu)]  \cdot  [(K, \psi), g_{2}, (H, \phi)] =   [(K, \psi), g_{1}g_{2}, (J, \mu)]  \]
and zero otherwise. This product makes sense because 

(i)  \   if $K \leq g_{2}^{-1} H g_{2}$ and
 $H \leq g_{1}^{-1} J g_{1}$ then $K \leq  g_{2}^{-1} H g_{2} \leq   g_{2}^{-1} g_{1}^{-1} J g_{1} g_{2} $ 
 
 and 
 
 (ii)  \  if $\psi(k) = \phi(h) = \mu(j), $ where $k = g_{2}^{-1}hg_{2},  h = g_{1}^{-1}j g_{1}$ then
 \linebreak
  $k = g_{2}^{-1}  g_{1}^{-1}j g_{1}   g_{2}$. 
 
 This product is clearly associative and we define an algebra ${\mathcal H}_{cmc}(G)$ to be ${\mathcal H}$ modulo the relations ( \cite{Sn18}\footnote{For the purposes of this essay I  am using throughout a convention in which $g$ is replaced by $g^{-1}$ to make the composition coincide with the conventions of induced representations which I shall use here.})
 \[   [(K, \psi), gk, (H, \phi)]  = \psi(k^{-1}) [(K, \psi), g, (H, \phi)]  \]
and
 \[     [(K, \psi), hg, (H, \phi)]  = \phi(h^{-1}) [(K, \psi), g, (H, \phi)] .    \]
  
 \section{The monocentre of a group}
 
 Suppose that we are in the situation of \S 2.
 
 To recapitulate:  Let $G$ be a locally profinite group and let $k$ be an algebraically closed field. Suppose that $\underline{\phi} : Z(G) \longrightarrow k^{*}$ is a fixed $k$-valued, continuous character on the centre $Z(G)$ of $G$. Let ${\mathcal M}_{cmc, \underline{\phi}}(G)$, as in \S2, be the poset of pairs $(H, \phi)$ where $H$ is a subgroup of $G$, containing $Z(G)$, which is compact, open modulo the centre of $G$ and $\phi : H \longrightarrow k^{*}$ is a $k$-valued, continuous character whose restriction to $Z(G)$ is $\underline{\phi}$.
 
 As $(K, \psi)$ varies over ${\mathcal M}_{cmc, \underline{\phi}}(G)$ suppose that we have a family of elements of $G$, $\{ x_{(K, \psi)}  \in  {\rm stab}_{G}(K, \psi) \}$ indexed by pairs $(K, \psi)$ where 
 $ {\rm stab}_{G}(K, \psi)$ denotes the stabiliser of  $(K, \psi) $  
\[   {\rm stab}_{G}(K, \psi) = \{ z \in G \ | \ zKz^{-1} = K , \psi(zkz^{-1}) = \psi(k) \ {\rm for \ all} \ k \in K \} .  \]This is equivalent to $K \leq x_{(K, \psi)}^{-1} K x_{(K, \psi)}$  and, for all $k \in K$,
\[    \psi( x_{(K, \psi)}^{-1} k x_{(K, \psi)}  )  = \psi(k) = x_{(K, \psi)}^{*}(\psi)( x_{(K, \psi)}^{-1} k x_{(K, \psi)}  ) \]
so that $[(K, \psi),  x_{(K, \psi)}  ,  (K, \psi)]$ is one of the basis vectors for ${\mathcal H}$ of \S 2.

Next suppose that $(H, \phi) \in {\mathcal M}_{cmc, \underline{\phi}}(G)$ and $x_{(H, \phi)}$ are similar data for another pair and that $ [(K, \psi), g   , (H, \phi)]$ is another basis element of ${\mathcal H}$ in \S 2.

The {\bf monocentre condition} relating these elements is defined by 
\vspace{4pt}

(i) \ $  g x_{(K, \psi)} g^{-1} \in {\rm stab}_{G}(H, \phi) $
\vspace{4pt}

and
\vspace{4pt}

(ii) \  $ g x_{(K, \psi)}g^{-1} = x_{(H, \phi)} \in  {\rm stab}_{G}(H, \phi)/{\rm Ker}(\phi) $.
\vspace{4pt}

 Observe that ${\rm Ker}(\phi)$ is a normal subgroup of ${\rm stab}_{G}(H, \phi)$. Therefore if 
 $ [(K, \psi), g   , (H, \phi)] $, $x_{(K, \psi)} $ and $x_{(H, \phi)} $ satisfy the monocentre condition then so do  $ [(K, \psi), g   , (H, \phi)] $, $x_{(K, \psi)}^{-1} $ and $x_{(H, \phi)}^{-1} $.
 
 Furthermore, if  $ [(K, \psi), g   , (H, \phi)] $, $x_{(K, \psi)} $ and $x_{(H, \phi)} $ satisfy the monocentre condition and $w \in {\rm Ker}(\psi) \leq K$ then $ [(K, \psi), g   , (H, \phi)] $, $x_{(K, \psi)}w $ and $x_{(H, \phi)}gwg^{-1} $ also satisfy the condition and $gwg^{-1}  \in {\rm Ker}(\phi) \leq H$.

 \begin{proposition}{$_{}$}
 \label{3.1}
 \begin{em}

 The monocentre condition implies that the two compositions, defined in \S2,
 \[   [(K, \psi), g   , (H, \phi)]  \cdot   [(K, \psi), x_{(K, \psi)}   , (K, \psi)]   \]
and
\[  [(H, \phi), x_{(H, \phi)}   , (H, \phi)] \cdot   [(K, \psi), g   , (H, \phi)]  \]
are equal in the algebra ${\mathcal H}_{cmc}(G)$ of \S 2.
 \end{em}
 \end{proposition}
 
{\bf Proof:}   \  By definition we have
\[   [(K, \psi), g   , (H, \phi)]  \cdot   [(K, \psi), x_{(K, \psi)}   , (K, \psi)]  =  [(K, \psi), g x_{(K, \psi)}  , (H, \phi)]  \]
and
\[  [(H, \phi), x_{(H, \phi)}   , (H, \phi)] \cdot   [(K, \psi), g   , (H, \phi)]  = [(K, \psi),  x_{(H, \phi)} g , (H, \phi)]  .\]

The monomial condition implies that there exists $z \in {\rm Ker}(\phi) \subseteq H$ such that
 \[  g^{-1} x_{(H, \phi)}^{-1}  =   x_{(K, \psi)}^{-1} g^{-1} z  \]
 or, equivalently, that
  \[  z  x_{(H, \phi)} g =    g x_{(K, \psi)} .   \]
By the relation defining ${\mathcal H}_{cmc}(G)$ at the end of \S 2
 \[  \begin{array}{ll}
   [(K, \psi), z x_{(H, \phi)} g , (H, \phi)]  &  = \phi(z)^{-1}     [(K, \psi), z x_{(H, \phi)} g , (H, \phi)] \\
   \\
   & =    [(K, \psi), z x_{(H, \phi)} g , (H, \phi)] . 
   \end{array}  \]
 $\Box$
 \begin{definition}
 \label{3.2}
 \begin{em}
  The monocentre of $G$, denoted by $Z_{{\mathcal M}}(G)$,  is the set of families 
  $\{ x_{(K, \psi)}  \in   {\rm stab}_{G}(K, \psi)/{\rm Ker}(\psi)\}$ such that for every $x_{(K, \psi)} $, $x_{(H, \phi)} $ and $g$ such that $   (K, \psi) \leq (g^{-1}Hg, (g)^{*}(\phi)) $ the monocentre condition holds, as introduced above.

  
 As we shall see in more detail, because the monocentre condition includes a central character which is common to the pairs $(K, \psi)$ and $(H, \phi)$, $Z_{{\mathcal M}}(G)$ is the product of subgroups 
 \end{em}
 \end{definition}
  \begin{remark}{The maximal pairs}
 \label{3.3}
 \begin{em}
 
 The elements in a family corresponding to an element of $Z_{{\mathcal M}}(G)$ are cosets such as $x_{(K, \psi)} \in {\rm stab}_{G}(K, \psi)/{\rm Ker}(\psi)$ as $(K, \psi)$ varies over ${\mathcal M}_{cmc, \underline{\phi}}(G)$.
 \end{em}
 \end{remark}
  If $(K , \psi)$ is maximal in the partially ordered set ${\mathcal M}_{cmc, \underline{\phi}}(G)$ and we have 
 \[    (K, \psi) \leq (g^{-1}Hg, (g)^{*}(\phi)) \]
 we must have $K = g^{-1}Hg$ and $\phi(gkg^{-1}) = \psi(k)$ for all $k \in K$. Then any pair
 \[ x_{(K, \psi)} \in {\rm stab}_{G}(K, \psi)/{\rm Ker}(\psi), \ x_{(H, \phi)} \in {\rm stab}_{G}(H, \phi)/{\rm Ker}(\phi)   \]
 in the family must satisfy $g x_{(K, \psi)} g^{-1} \in {\rm stab}_{G}(H, \phi)$ automatically and also that 
 \[  g x_{(K, \psi)}g^{-1} = x_{(H, \phi)} \in  {\rm stab}_{G}(H, \phi)/{\rm Ker}(\phi). \]
 This means that the elements in a monocentre family corresponding to maximal pairs are determined by any one maximal pair - the other entries being conjugates as shown above. - which depend only on $g$'s modulo the stabiliser, which is precisely the indeterminacy which the conditions impose on $g$.
 
 Notwithstanding the above remark, it remains to decide whether or not $Z_{{\mathcal M}}(G)$ can be non-trivial.

 \begin{example}{$G = D_{8}$}
\label{3.4}
\end{example} 
The dihedral group of order eight is given by
\[  D_{8} = \langle x, y \ | \ x^{4} = 1 = y^{2}, \ yxy = x^{3} \rangle .  \]

The subgroups and the one-dimensional complex characters on them are
\[ \begin{array}{l}
D_{8}: \ 1, \chi_{1}, \chi_{2}, \chi_{1} \chi_{2}  \\
\hspace{40pt}  \chi_{1}(x) = -1, \chi_{1}(y)=1, \chi_{2}(x)=1, \chi_{2}(y) = -1 \\
\langle x \rangle: \  1, \phi, \phi^{2}, \phi^{3} \\
\hspace{40pt}  \phi(x) = i = \sqrt{-1}\\
\langle x^{2}, y \rangle:  \  1, \tilde{\chi}_{1}, \tilde{\chi}_{2}, \tilde{\chi}_{1} \tilde{\chi}_{2}  \\
\hspace{40pt}  \tilde{\chi}_{1}(x^{2}) = -1, \tilde{\chi}_{1}(y)=1, \tilde{\chi}_{2}(x^{2})=1, \tilde{\chi}_{2}(y) = -1 \\
\langle x^{2}, xy \rangle : \  1, \hat{\chi}_{1}, \hat{\chi}_{2}, \hat{\chi}_{1} \hat{\chi}_{2}  \\
\hspace{40pt}  \hat{\chi}_{1}(x^{2}) = -1, \hat{\chi}_{1}(xy)=1, \hat{\chi}_{2}(x^{2})=1, \hat{\chi}_{2}(xy) = -1 \\
\langle x^{2} \rangle : \  1, \chi  \\
\langle 1 \rangle : \  1 \\
\end{array} \]

Therefore the elements of ${\mathcal M}_{cmc, 1}(D_{8})$ are
\[ \begin{array}{l}
D_{8}: \ 1, \chi_{1}, \chi_{2}, \chi_{1} \chi_{2}  \\
\langle x \rangle: \  1,  \phi^{2},  \\
\langle x^{2}, y \rangle:  \  1, \tilde{\chi}_{2}  \\
\langle x^{2}, xy \rangle : \  1,  \hat{\chi}_{2}  \\
\langle x^{2} \rangle : \  1  \\
\end{array} \]

 Now consider  
 \[ x_{(\langle x^{2} \rangle , 1)}  \in {\rm stab}_{D_{8}}(\langle x^{2} \rangle,  1)/{\rm Ker}( 1 ) = D_{8}/\langle x^{2} \rangle.\]
We have inclusions
\[ \begin{array}{l}
 (\langle x^{2} \rangle , 1) \leq  (D_{8}, 1) \\
 (\langle x^{2} \rangle , 1) \leq (D_{8}, \chi_{1}) \\
 (\langle x^{2} \rangle , 1) \leq  (D_{8}, \chi_{2}) \\
 (\langle x^{2} \rangle , 1) \leq (D_{8}, \chi_{1} \chi_{2})  \\
(\langle x^{2} \rangle , 1) \leq  \langle x \rangle ,  1)  \\
(\langle x^{2} \rangle , 1) \leq ( \langle x \rangle,  \phi^{2})  \\
(\langle x^{2} \rangle , 1) \leq  (\langle x^{2}, y \rangle, 1)  \\
(\langle x^{2} \rangle , 1) \leq   (\langle x^{2}, y \rangle, \tilde{\chi}_{2} ) \\
(\langle x^{2} \rangle , 1) \leq  (\langle x^{2}, xy \rangle ,  1)  \\
(\langle x^{2} \rangle , 1) \leq  ,  (\langle x^{2}, xy \rangle ,  \hat{\chi}_{2})  \\
\end{array} \]
which require monocentre conditions which are depicted as follows:
\[ \begin{array}{l}
x_{ (\langle x^{2} \rangle , 1)} \mapsto x_{ (D_{8}, 1)} ,   \ D_{8}/\langle x^{2} \rangle \longrightarrow \{ 1 \}, \\
 x_{(\langle x^{2} \rangle , 1)} \mapsto  x_{ (D_{8}, \chi_{1})}, \  D_{8}/\langle x^{2} \rangle   \longrightarrow D_{8}/\langle x^{2}, y \rangle  , \\
 x_{(\langle x^{2} \rangle , 1)} \mapsto x_{ (D_{8}, \chi_{2})},   \ D_{8}/\langle x^{2} \rangle \longrightarrow D_{8}/\langle x \rangle , \\
 x_{(\langle x^{2} \rangle , 1)} \mapsto x_{ (D_{8}, \chi_{1} \chi_{2})},   \ D_{8}/\langle x^{2} \rangle \longrightarrow  D_{8}/\langle x^{2}, xy \rangle , \\
 x_{(\langle x^{2} \rangle , 1)} \mapsto x_{  \langle x \rangle ,  1)}   \ D_{8}/\langle x^{2} \rangle \longrightarrow D_{8}/\langle x \rangle , \\
 x_{(\langle x^{2} \rangle , 1)} \mapsto x_{( \langle x \rangle,  \phi^{2})}  ,  \ D_{8}/\langle x^{2} \rangle \longrightarrow  D_{8}/\langle x^{2} \rangle , \\
 x_{(\langle x^{2} \rangle , 1)} \mapsto x_{(\langle x^{2}, y \rangle,  1) },  \ D_{8}/\langle x^{2} \rangle \longrightarrow D_{8}/\langle x^{2}, y \rangle , \\
 x_{(\langle x^{2} \rangle , 1)} \mapsto  x_{(\langle x^{2}, y \rangle, \tilde{\chi}_{2} )},  \ D_{8}/\langle x^{2} \rangle \longrightarrow  D_{8}/\langle x^{2} \rangle , \\
 x_{(\langle x^{2} \rangle , 1)} \mapsto  x_{ (\langle x^{2}, xy \rangle ,  1)},  \ D_{8}/\langle x^{2} \rangle \longrightarrow D_{8}/\langle x^{2}, xy \rangle ,  \\
 x_{(\langle x^{2} \rangle , 1)} \mapsto x_{(\langle x^{2}, xy \rangle ,  \hat{\chi}_{2})} ,  \ D_{8}/\langle x^{2} \rangle \longrightarrow  D_{8}/\langle x^{2} \rangle .\\
\end{array} \]

Therefore $Z_{ {\mathcal M}_{cmc, 1}}(D_{8}) \cong  D_{8}/\langle x^{2} \rangle$.

Next consider   
 \[ x_{(\langle x^{2} \rangle , \chi)}  \in {\rm stab}_{D_{8}}(\langle x^{2} \rangle,  \chi)/{\rm Ker}( \chi ) = D_{8}.\]

The elements of ${\mathcal M}_{cmc, \chi}(D_{8})$ are
\[ \begin{array}{l}
\langle x \rangle: \   \phi,  \phi^{3} \\
\langle x^{2}, y \rangle:  \  \tilde{\chi}_{1}, \tilde{\chi}_{1} \tilde{\chi}_{2}  \\
\langle x^{2}, xy \rangle : \  \hat{\chi}_{1},  \hat{\chi}_{1} \hat{\chi}_{2}  \\
\langle x^{2} \rangle : \   \chi  \\
\end{array} \]
so we have inclusions
\[  \begin{array}{l}
(\langle x^{2} \rangle , \chi) \leq  ( \langle x \rangle,  \phi) , \\
(\langle x^{2} \rangle , \chi) \leq (\langle x \rangle, \phi^{3}) \\
(\langle x^{2} \rangle , \chi) \leq  ( \langle x^{2}, y \rangle,  \tilde{\chi}_{1}),   \\
(\langle x^{2} \rangle , \chi) \leq   (\langle x^{2}, y \rangle , \tilde{\chi}_{1} \tilde{\chi}_{2} ), \\
(\langle x^{2} \rangle , \chi) \leq ( \langle x^{2}, xy \rangle ,  \hat{\chi}_{1}),  \\
(\langle x^{2} \rangle , \chi) \leq  (\langle x^{2}, xy \rangle,  \hat{\chi}_{1} \hat{\chi}_{2}).  \\
\end{array} \]
which require monocentre conditions which are depicted as follows:
\[ \begin{array}{l}
x_{(\langle x^{2} \rangle , \chi)} \longrightarrow  x_{ ( \langle x \rangle,  \phi)} ,  \  D_{8} \longrightarrow  \langle x \rangle \\
x_{(\langle x^{2} \rangle , \chi)} \longrightarrow  x_{(\langle x \rangle, \phi^{3})},  \ D_{8} \longrightarrow   \langle x \rangle  \\
x_{(\langle x^{2} \rangle , \chi)} \longrightarrow  x_{ ( \langle x^{2}, y \rangle,  \tilde{\chi}_{1})}, \ D_{8} \longrightarrow   \langle x^{2}, y \rangle/   \langle  y \rangle  \\
x_{(\langle x^{2} \rangle , \chi)} \longrightarrow  x_{   (\langle x^{2}, y \rangle , \tilde{\chi}_{1} \tilde{\chi}_{2} )}, \ D_{8} \longrightarrow   \langle x^{2}, y \rangle/ \langle x^{2} y \rangle  \\
x_{(\langle x^{2} \rangle , \chi)} \longrightarrow  x_{( \langle x^{2}, xy \rangle ,  \hat{\chi}_{1})},  \ D_{8} \longrightarrow   \langle x^{2}, xy \rangle/  \langle xy \rangle   \\
x_{(\langle x^{2} \rangle , \chi)}    \longrightarrow  x_{  (\angle x^{2}, xy \rangle,  \hat{\chi}_{1} \hat{\chi}_{2})}, \ D_{8} \longrightarrow  \langle x^{2}, xy \rangle/  \langle x^{3} y \rangle   .  \\
\end{array} \]
In the case of the first two conditions in the list we see that 
\[ x_{(\langle x^{2} \rangle , \chi)} \in \langle x \rangle \]
is necessary. However, in order to satisfy the monocentre condition in relation to the last four entries in the list we require that
\[ x_{(\langle x^{2} \rangle , \chi)} \in \langle x^{2} \rangle .\]
Therefore we obtain 
\[ Z_{ {\mathcal M}}(D_{8}) = 
Z_{ {\mathcal M}_{cmc, 1}}(D_{8})  \times Z_{ {\mathcal M}_{cmc, \chi}}(D_{8}) \cong 
D_{8}/ \langle x^{2} \rangle   \times   \langle x^{2} \rangle. \]

This dihedral example reveals the following result.
\begin{theorem}{$_{}$}
\label{3.5}
\begin{em}

The monocentre group, $Z_{{\mathcal M}}(G)$, is the product of the subgroups $Z_{ {\mathcal M}_{cmc, \underline{\phi}}}(G)$ as $\underline{\phi}$ varies over the central characters. Also the set of elements in a family   $\{ x_{(K, \psi)}  \in   {\rm stab}_{G}(K, \psi)/{\rm Ker}(\psi)\}$ representing an element of 
$Z_{ {\mathcal M}_{cmc, \underline{\phi}}}(G)$ are determined by the 
\[    x_{(Z(G), \underline{\phi})}  \in   G/{\rm Ker}(\underline{\phi})   \]
such that the image of $x_{(Z(G), \underline{\phi})} $ represents an element 
\linebreak
$x_{(K, \psi)}  \in  {\rm stab}_{G}(K, \psi)/{\rm Ker}(\psi)$ for every $(K, \phi) \in 
 {\mathcal M}_{cmc, \underline{\phi}}$.
\end{em}
\end{theorem}
\begin{remark}
\label{3.6}
\begin{em}

(i) \ The monocentre group is an entertaining construction, but it will turn out to be too restrictive for our purposes. Although it might be less trivial - even useful! - in the case of modular representations.

(ii) \ More important is the following situation (see \S\ref{11.1.3} and \S\ref{11.1.4}). Fix a central character $\underline{\phi}$ as usual.
\end{em}
\end{remark}

Suppose, for $i=1,2$, that we are given 
\[ \begin{array}{l}
[(K_{i}, \psi_{i}), g_{i}, (H_{i}, \phi_{i})] \ {\rm   and}  \\
\\
 \{ x_{(K_{i}, \psi_{i})}  \in   {\rm stab}_{G}(K_{i}, \psi_{i})/{\rm Ker}(\psi_{i})  \}   \ {\rm   and}  \\
 \\
  \{ x_{(H_{i}, \phi_{i})}  \in   {\rm stab}_{G}(H_{i}, \phi_{i})/{\rm Ker}(\phi_{i})  \}  
\end{array} \]
which satisfy both
\[   \begin{array}{l}
 [(H_{1}, \phi_{1}), x_{(H_{1}, \phi_{1})}, (H_{1}, \phi_{1})] \cdot 
 [(K_{1}, \psi_{1}), g_{1}, (H_{1}, \phi_{1})]   \\
 \\
 =    
   [(K_{1}, \psi_{1}), g_{1}, (H_{1}, \phi_{1})]    \cdot   [(K_{1}, \psi_{1}), x_{(K_{1}, \psi_{1})}, (K_{1}, \psi_{1})] 
     \end{array} \]
     and
  \[   \begin{array}{l}
 [(H_{2}, \phi_{2}), x_{(H_{2}, \phi_{2})}, (H_{2}, \phi_{2})] \cdot 
 [(K_{2}, \psi_{2}), g_{2}, (H_{2}, \phi_{2})]   \\
 \\
 =    
   [(K_{2}, \psi_{2}), g_{2}, (H_{2}, \phi_{2})]    \cdot  
    [(K_{2}, \psi_{2}), x_{(K_{2}, \psi_{2})}, (K_{2}, \psi_{2})] .
     \end{array} \]   
     
 Under these conditions we require that for all 
 \[    [(H_{1}, \phi_{1}), g_{3}, (H_{2}, \phi_{2})]  \ {\rm and} \    [(K_{1}, \psi_{1}), g_{4}, (K_{2}, \psi_{2})]  \]
 such that
 \[  \begin{array}{l}
   [(H_{1}, \phi_{1}), g_{3}, (H_{2}, \phi_{2})]  \cdot [(K_{1}, \psi_{1}), g_{1}, (H_{1}, \phi_{1})] \\
   \\
    =
  [(K_{2}, \psi_{2}), g_{2}, (H_{2}, \phi_{2})]   \cdot  [(K_{1}, \psi_{1}), g_{4}, (K_{2}, \psi_{2})] 
  \end{array}  \]
  the $\{ x_{(K_{i}, \psi_{i})} , x_{(H_{i}, \phi_{i})} \}$ satisfy
\[  \begin{array}{l}  
 [(H_{2}, \phi_{2}), x_{(H_{2}, \phi_{2})}, (H_{2}, \phi_{2})]  \cdot    [(H_{1}, \phi_{1}), g_{3}, (H_{2}, \phi_{2})]   \\
 \\
 =    [(H_{1}, \phi_{1}), g_{3}, (H_{2}, \phi_{2})]   \cdot  [(H_{1}, \phi_{1}), x_{(H_{1}, \phi_{1})}, (H_{1}, \phi_{1})]  \\
   \end{array}  \]
and also that 
 \[ \begin{array}{l}
   [(K_{2}, \psi_{2}), x_{(K_{2}, \psi_{2})}, (K_{2}, \psi_{2})]  \cdot   [(K_{1}, \psi_{1}), g_{4}, (K_{2}, \psi_{2})]   \\
  \\
=    [(K_{1}, \psi_{1}), g_{4}, (K_{2}, \psi_{2})] 
    \cdot   [(K_{1}, \psi_{1}), x_{(K_{1}, \psi_{1})}, (K_{1}, \psi_{1})] .
  \end{array}  \]

\section{Extending the definition of admissibility}

If $G$ is a locally profinite group and $k$ is an algebraically closed field then a $k$-representation of $G$ is a vector space $V$ with a left, $k$-linear $G$-action. Let $\underline{\phi} : Z(G) \longrightarrow k^{*}$ be a continuous character on the centre of $G$. Let ${\mathcal M}_{cmc, \underline{\phi}}(G)$, as in \S2, denote the poset of pairs $(H, \phi)$ where $H$ is a subgroup of $G$, such that $Z(G) \subseteq H$, which is compact open modulo the centre and $\phi : H \longrightarrow k^{*}$ is a continuous character which extends $\underline{\phi}$.

Suppose that $V$ is acted upon by $g \in Z(G)$ via multiplication by $\underline{\phi}(g)$. The representation $V$ is called smooth if
\[ V = \bigcup_{K \subset G,  \ K \ {\rm compact, open}} \ V^{K} .\]
$V$ is called admissible if ${\dim}_{k}( V^{K}) < \infty$ for all compact open subgroups $K$. Define a subspace of $V$, denoted by $V^{(H, \phi)}$, for $(H, \phi) \in {\mathcal M}_{cmc, \underline{\phi}}(G)$ by
 \[  V^{(H, \phi)} = \{ v \in V \ | \ g \cdot v =  \phi(g)v \ {\rm for \ all} \ g \in H \}   . \]
 Hence $V^{K} = V^{(Z(G) \cdot K,  \phi)}$ if $\phi$ is a continuous character which is trivial on $K$.
 
 We shall say that $V$ is ${\mathcal M}_{cmc, \underline{\phi}}(G)$-smooth if 
 \[ V = \bigcup_{ (H, \phi) \in {\mathcal M}_{cmc, \underline{\phi}}(G)} \ V^{(H, \phi)}.\]
 In addition we shall say that $V$ is ${\mathcal M}_{cmc, \underline{\phi}}(G)$-admissible if 
 ${\rm dim}_{k}V^{(H, \phi)} < \infty$ for all  $(H, \phi) \in {\mathcal M}_{cmc, \underline{\phi}}(G)$.

 \begin{proposition}{$_{}$}
 \label{4.1}
 \begin{em}
 
 Let $G$ be a locally profinite group and let $k$ be an algebraically closed field. Let $V$ be a $k$-representation of $G$ with central character $\underline{\phi}$. Suppose that every continuous, $k$-valued character of a compact open subgroup of $G$ has finite image. Then $V$ is 
 ${\mathcal M}_{cmc, \underline{\phi}}(G)$-admissible if and only if it is admissible.
 \end{em}
 \end{proposition}
 \vspace{5pt}
 
 {\bf Proof}:
  \vspace{5pt}
  
  If $K$ is compact open then $K \bigcap Z(G)$ is also compact open. It is certainly compact, being a closed subset of a compact subspace. For $G = GL_{n}F$ with $F$ a $p$-adic local field the assumption it true. More generally, it holds if the quotient of $Z(G)$ by its maximal compact subgroup is 
discrete\footnote{If this condition is not true in general it is true in the main cases of interest. Therefore let us treat it as an unimportant assumption for the time being!}.
  
 Suppose that $V$ is admissible. If $H$ is a subgroup of $G$ which is compact open modulo the centre then $H = Z(G) \cdot K$ for some compact open subgroup. In this case supose that $\phi$ is a character of $H$ extending the central character. Then $V^{(H, \phi)}  = V^{(K, \mu)}$ where $\mu = {\rm Res}_{K}^{H}(\phi)$. Since the image of $\mu$ is finite the kernel of $\mu$ is compact open and 
 $V^{(K, \mu)} \subseteq V^{{\rm Ker}(\mu)}$, which is finite-dimensional.
 
 Next suppose that $0 \not= v \in V$. There exists a compact open subgroup $K$ such that $v \in V^{K}$.
 Set $H = Z(G) \cdot K$, which is compact open modulo $Z(G) \subset H$. If $g \in Z(G) \bigcap K$ then
 $v = g \cdot v = \underline{\phi}(g) \cdot v$ so that the central character is trivial on $Z(G) \bigcap K$. 
 Hence the central character induces a character $\lambda$ on $H$ which factors through $K/Z(G) \bigcap K \cong Z(G) \cdot K/K$ and so $v \in V^{(H, \lambda)}$, which completes the proof of 
 ${\mathcal M}_{cmc, \underline{\phi}}(G)$-admissibility.
 
 Assume that $V$ is ${\mathcal M}_{cmc, \underline{\phi}}(G)$-admissible. If $0 \not= v \in V$ belongs to 
 $V^{(H, \phi)}$ where $H$ is compact open modulo the centre then $H = Z(G) \cdot K$ where $K$ is compact open. Hence $v \in V^{J}$ where $J$ is the compact open subgroup given by 
 $J = {\rm Ker}( {\rm Res}_{K}^{H}(\phi))$. 
 
 Next suppose that $K$ is a compact open subgroup. If $V^{K}$ is non-trivial then $V^{K} \subseteq V^{(Z(G) \cdot K, \lambda)}$ where $\lambda : H = Z(G) \cdot K \longrightarrow k^{*}$ is the character which was constructed in the first half of the proof. Since  $V^{(Z(G) \cdot K, \lambda)}$ is assumed to be finite-dimensional this concludes the proof of admissibility. $\Box$
\begin{question}{Di-$p$-adic Langlands}
\label{4.2}
\begin{em}

In the last 20 years I believe that several authors have studied the ``$p$-adic Langlands programme''. This is the situation where, for example, one studies ``admissible'' representations of a locally $p$-adic Lie group on vector spaces over the algebraic closure of a $p$-adic local field (or its residue field).

I intend to call this the di-$p$-adic situation since it is no more complicated to say and indicates the involvement of $p$-adic fields twice. In addition to \cite{AKdS13} there are lots of papers on this subject\footnote{Regrettably I have not got round to reading any of them!} and a useful source for these (brought to my attention by Rob Kurinczuk) is the bibliography of  \cite{CEGGPS14}.

The question arises: Are the sort of representations considered by the di-$p$-adic professionals 
${\mathcal M}_{cmc, \underline{\phi}}(G)$-admissible?

The first issue, which I expect has been studied a lot already in the literature, is the ``action'' of the finite-dimensional di-$p$-adic representations on the infinite-dimensional ones. 

I guess that the classification of the finite-dimensional di-$p$-adic representations for general linear groups follows the classical pattern?
\end{em}
\end{question}

\section{Induced representations and locally profinite groups}

Let $G$ be a locally profinite group. In this section we are going to study admissible representations of $G$ and its subgroups in relation to induction. These representations will be given by left-actions of the groups on vector spaces over $k$, which is an algebraically closed field of arbitrary characteristic. 

Let us begin by recalling, from (\cite{Sn18} Chapter Two \S1), induced and compactly induced smooth representations.

\begin{definition}{$_{}$}
\label{5.1}
\begin{em}

Let $G$ be a locally profinite group and $H \subseteq G$ a closed subgroup. Thus $H$ is also  locally profinite. Let 
\[   \sigma : H \longrightarrow  {\rm Aut}_{k}(W)    \]
be a smooth representation of $H$. Set $X$ equal to the space of functions $f: G \longrightarrow W  $ such that (writing simply $h \cdot w$ for $\sigma(h)(w)$ if $h \in H, w \in W$)

(i) \  $f(hg) = h \cdot f(g)$ for all $h \in H, g \in G$,

(ii)  \  there is a compact open subgroup $K_{f} \subseteq G$ such that $f(gk) = f(g)$ for all $g \in G, k \in K_{f}$.

The (left) action of $G$ on $X$ is given by $(g \cdot f)(x)= f(xg)$ and
\[   \Sigma :  G  \longrightarrow  {\rm Aut}_{k}(X)  \]
gives a smooth representation of $G$.

The representation $\Sigma$ is called the representation of $G$ smoothly induced from $\sigma$ and is usually denoted by $\Sigma = {\rm Ind}_{H}^{G}(\sigma)$.
\end{em}
\end{definition}

\begin{dummy}
\label{5.2}
\begin{em}

\[  (g \cdot f)(hg_{1}) = f(hg_{1}g) = h f(g_{1}g) = h (g \cdot f)(g_{1}) \]  
so that $(g \cdot f)$ satisfies condition (i) of Definition \ref{5.1}. 

Also, by the same discussion as in the finite group case (see Appndix \S4), the formula will give a left $G$-representation, providing that $g \cdot f \in X$ when $f \in X$.
However, condition (ii) asserts that there exists a compact open subgroup $K_{f}$ such that $k \cdot f = f$ for all $k \in K_{f}$. The subgroup $gK_{f}g^{-1}$ is also a compact open subgroup and, if $k \in K_{f}$,
we have
\[ (gkg^{-1}) \cdot (g \cdot f) = (gkg^{-1}g) \cdot f =  (gk) \cdot f = (g \cdot (k \cdot f)) =  (g \cdot f) \]
so that $g \cdot f \in X$, as required.

The smooth representations of $G$ form an abelian category ${\rm Rep}(G)$.
\end{em}
\end{dummy}

\begin{proposition}{$_{}$}
\label{5.3}
\begin{em}

The functor
\[  {\rm Ind}_{H}^{G} : {\rm Rep}(H) \longrightarrow   {\rm Rep}(G)  \]
is additive and exact.
\end{em}
\end{proposition}
\begin{proposition}{(Frobenius Reciprocity) }
\label{5.4}
\begin{em}

There is an isomorphism
\[  {\rm Hom}_{G}( \pi,  {\rm Ind}_{H}^{G}(\sigma)) \stackrel{\cong}{\longrightarrow}  {\rm Hom}_{H}( \pi,  \sigma)  \]
given by $\phi \mapsto \alpha \cdot \phi$ where $\alpha$ is the $H$-map  
\[   {\rm Ind}_{H}^{G}(\sigma)  \longrightarrow  \sigma  \]
given by $\alpha(f) = f(1)$.
\end{em}
\end{proposition}

\begin{dummy}
\label{5.5}
\begin{em}

In general, if $H \subseteq Q$ are two closed subgroups there is a $Q$-map
\[   {\rm Ind}_{H}^{G}(\sigma)  \longrightarrow   {\rm Ind}_{H}^{Q}(\sigma)   \]
given by restriction of functions. Note that $\alpha$ in Proposition \ref{5.4} is the special case where $H=Q$.
\end{em}
\end{dummy}

\begin{dummy}{The c-Ind variation}
\label{5.6}
\begin{em}

Inside $X$ let $X_{c}$ denote the set of functions which are compactly supported modulo $H$. This means that the image of the support
\[   {\rm supp}(f) = \{ g \in G \ | \  f(g) \not= 0  \}  \]
has compact image in $H \backslash G$. Alternatively  there is a compact subset $C \subseteq G$ such that $  {\rm supp}(f) \subseteq H \cdot C$.

The $\Sigma$-action on $X$ preserves $X_{c}$, since $ {\rm supp}(g \cdot f) =  {\rm supp}(f) g^{-1} \subseteq HCg^{-1}$, and we obtain $X_{c} =  c- {\rm Ind}_{H}^{G}(W) $, the compact induction of $W$ from $H$ to $G$.

This construction is of particular interest when $H$ is open. There is a canonical left $H$-map (see the Appendix  in induction in the case of finite groups)
\[    f : W \longrightarrow  c- {\rm Ind}_{H}^{G}(W)  \]
given by $w \mapsto f_{w}$ where $f_{w}$ is supported in $H$ and $f_{w}(h) = h \cdot w$ (so $f_{w}(g) =0$ if $g \not\in H$).

For $g \in G$ we have
\[ \begin{array}{ll}
(g \cdot f_{w})(x) = f_{w}(xg) & =  \left\{  \begin{array}{cc}
0 & {\rm if } \ xg \not\in H, \\
\\
(xg^{-1}) \cdot w & {\rm if } \ xg \in H, \\
\end{array} \right.   \\
\\
& =  \left\{  \begin{array}{cc}
0 & {\rm if } \  x  \not\in Hg^{-1}, \\
\\
(xg^{-1}) \cdot w & {\rm if } \ x \in Hg^{-1}. \\
\end{array} \right. 
\end{array} \]

We shall be particularly interested in the case when ${\rm dim}_{k}(W) =1$. In this case we write $W = k_{\phi}$ where $\phi : H \longrightarrow k^{*}$ is a continuous/smooth character and, as a vector space with a left $H$-action $W =k$ on which $h \in H$ acts by multiplication by $\phi(h)$. In this case $\alpha_{c}$ is an
injective left $k[H]$-module homomorphism of the form
\[   f : k_{\phi} \longrightarrow   c- {\rm Ind}_{H}^{G}(k_{\phi}) .\]

\end{em}
\end{dummy}

\begin{lemma}{$_{}$}
\label{5.7}
\begin{em}

Let $H$ be an open subgroup of $G$. Then

(i) \  $f : w \mapsto f_{w}$ is an $H$-isomorphism onto the space of functions $f \in c- {\rm Ind}_{H}^{G}(W) $ such that $  {\rm supp}(f)  \subseteq H$.

(ii) \  If $w \in W$ and $h \in H$ then $   h \cdot f_{w}   =  f_{h \cdot w}$.

(iii)  \   If ${\mathcal W}$ is a $k$-basis of $W$ and ${\mathcal G}$ is a set of coset representatives for 
$H \backslash G  $ then
\[   \{  g  \cdot f_{w}  \  |  \  w \in {\mathcal W}, \   g \in {\mathcal G}  \}\]
is a $k$-basis of $c- {\rm Ind}_{H}^{G}(W) $.
\end{em}
\end{lemma}
\vspace{2pt}

{\bf Proof}
\vspace{2pt}

If ${\rm supp}(f)$ is compact modulo $H$ there exists a compact subset $C$ such that 
\[ {\rm supp}(f) \subseteq HC  = \bigcup_{c \in C} \ Hc. \]
Each $Hc$ is open  so the open covering of $C$ by the $Hc$'s refines to a finite covering and so
\[   C  =  Hc_{1} \bigcup \ldots \bigcup Hc_{n}  \]
and so
\[ {\rm supp}(f) \subseteq HC  = Hc_{1} \bigcup \ldots \bigcup Hc_{n} . \]

For part (i), the map $f$ is an $H$-homomorphism to the space of functions supported in $H$ with inverse map $f \mapsto f(1)$.

For part (ii), from \S\ref{2.1.6} we have
\[  (h \cdot f_{w} )(x) =  f_{w}(xh) =  \left\{  \begin{array}{ll}
0 & {\rm if} \ x   \not\in H, \\
\\
xh \cdot  w & {\rm if} \   x  \in H. \\
\end{array} \right.   \]
so that, for all $x \in G$, $(h \cdot f_{w} )(x) = f_{h \cdot  w }(x)$, as required.

For part (iii), the support of any $f  \in c- {\rm Ind}_{H}^{G}(W) $ is a finite union of cosets $Hg$ where the $g$'s are chosen from the set of coset representatives ${\mathcal G}$ of $H \backslash G $.
The restriction of $f$ to any one of these $Hg$'s also lies in $c- {\rm Ind}_{H}^{G}(W) $. If ${\rm supp}(f) \subseteq Hg$ then $(g \cdot f)(z) \not= 0$ implies that $zg \in Hg$ so that $g \cdot f$ has support contained in $H$. Hence $g \cdot f$ on $H$ is a finite linear combination of the functions $f_{w}$ with $w \in {\mathcal W}$. Therefore $f$ is a finite linear combination of 
$g  \cdot f_{w}$'s where $  w \in {\mathcal W},    g \in {\mathcal G} $.
Clearly the set of functions $g \cdot  f_{w}$ with $g \in {\mathcal G}$ and $w \in {\mathcal W}$ is linearly independent. $\Box$ 

\begin{example}
\label{5.8}
\begin{em}

Let $K$ be a $p$-adic local field with valuation ring ${\mathcal O}_{K}$ and $\pi_{K}$ a generator of the maximal ideal of ${\mathcal O}_{K}$. Suppose that $G = GL_{n}K$ and that $H$ is a subgroup containing the centre of $G$ (that is, the scalar matrices $K^{*}$).  If $H$ is compact, open modulo $K^{*}$ then there is a subgroup $H'$ of finite index in $H$ such that $H' = K^{*}H_{1}$ with $H_{1}$ compact, open in $SL_{n}K$. This can be established by studying the simplicial action of $GL_{n}K$ on a suitable barycentric subdivision of the Bruhat-Tits building of $SL_{n}K$ (see \cite{Sn18} Chapter Four \S1).

To show that $H$ is both open and closed it suffices to verify this for $H'$. Firstly $H'$ is open, since it is $H' = \bigcup_{z \in K^{*}} \ zH_{1} = \bigcup_{s \in {\mathbb Z}}  \  \pi_{K}^{s} H_{1}$. 

Also $H' = K^{*}H_{1}$ is closed. Suppose that $X'  \not\in   K^{*}H_{1}$. $K^{*}H_{1}$ is closed under mutiplication by the multiplicative group generated by $\pi_{K}$ so that $\pi_{K}^{m} X' \not\in  K^{*}H_{1}$ for all $m$. By conjugation we may assume that $H_{1}$ is a subgroup of $SL_{n}{\mathcal O}_{K}$, which is the maximal compact open subgroup of $SL_{n}K$, unique up to conjugacy. Choose the smallest non-negative integer $m$ such that every entry of $X = \pi_{K}^{m} X' $ lies in
$ {\mathcal O}_{K}$. Therefore we may write $0 \not= {\rm det}(X) = \pi_{K}^{s} u$ where $u \in {\mathcal O}_{K}^{*}$ and $1 \leq s$. Now suppose that $V$ is an $n \times  n$ matrix with entries in $ {\mathcal O}_{K}$ such that $X + \pi_{K}^{t}V  \in K^{*}H_{1}$. Then
\[ {\rm det}(X + \pi_{K}^{t}V)  \equiv  \pi_{K}^{s} u  \ ({\rm modulo} \  \pi_{K}^{t}) . \]
So that if $t > s$ then $s$ must have the form $s = nw$ for some integer $w$ and 
$\pi_{K}^{-w} (X + \pi_{K}^{t}V) \in GL_{n}{\mathcal O}_{K} \bigcap K^{*}H_{1} = H_{1}$. Therefore all the entries in $\pi_{K}^{-w} X$ lie in ${\mathcal O}_{K} $ and $\pi_{K}^{-w} X \in GL_{n}{\mathcal O}_{K}$. Enlarging $t$, if necessary, we can ensure that $\pi_{K}^{-w} X \in H_{1}$, since $H_{1}$ is closed (being compact), and therefore 
$X' \in K^{*} H_{1}$, which is a contradiction.

Since $H $ is both closed and open in $GL_{n}K$ we may form the admissible representation $c-{\rm Ind}_{H}^{GL_{n}K}(k_{\phi})$ for any continuous character $\phi : H \longrightarrow k^{*}$ and apply Lemma \ref{2.1.7}. 

If $g \in GL_{n}K, h \in H$ then $(g \cdot f_{1})(x) =  \phi(xg)$ if $xg \in H$ and zero otherwise. On the other hand, $(gh \cdot f_{1})(x) =  \phi(xgh) = \phi(h)\phi(xg)$ if $xg \in H$ and zero otherwise. Therefore as a left $GL_{n}K$-representation $c-{\rm Ind}_{H}^{GL_{n}K}(k_{\phi})$
is isomorphic to 
\[ k[GL_{n}K] /( \phi(h) g  -  gh \ | \  g \in GL_{n}K, \ h \in H)  \]
with left action induced by $g_{1} \cdot g = g_{1}g$. 

This vector space is isomorphic to the $k$-vector space whose basis is given by $k$-bilinear tensors over $H$ of the form $g \otimes_{k[H]} 1$ as in the case of finite groups. The basis vector $g \cdot f_{1}$ corresponds to $g \otimes_{H} 1$ and $GL_{n}K$ acts on the tensors by left multiplication, as usual (see Appendix \S4 in the finite group case).
\end{em}
\end{example}

\begin{proposition}{$_{}$}
\label{5.9}
\begin{em}

The functor
\[   c- {\rm Ind}_{H}^{G} : {\rm Rep}(H) \longrightarrow   {\rm Rep}(G)  \]
is additive and exact.
\end{em}
\end{proposition}
\newpage

\begin{proposition}{$_{}$}
\label{5.10}
\begin{em}

Let $H \subseteq G$ be an open subgroup and $(\sigma ,  W)$ smooth. Then there is a functorial isomorphism
\[  {\rm Hom}_{G}(  c-{\rm Ind}_{H}^{G}(W) , \pi) \stackrel{\cong}{\longrightarrow}  
{\rm Hom}_{H}(  W  , \pi )  \]
given by $F \mapsto F \cdot  f $, the composition with the $H$-map $f$ of Lemma \ref{5.7}.
\end{em}
\end{proposition}

\begin{example}{$c-\underline{{\rm Ind}}_{H}^{G}(\phi)$}
\label{5.11}
\begin{em}

Suppose that $\phi : H \longrightarrow k^{*}$ is a continuous character (i.e. a one-dimensional smooth representation of $H$).

Suppose that we are in a situation analogous to that of Example \ref{5.8}. Namely suppose that $H$ is open and closed, contains $Z(G)$,the centre of $G$, and is compact open modulo $Z(G)$. A basis for $k$ is given by $1 \in k^{*}$ and we have the function $f_{1} \in X_{c}$ given by
$f_{1}(h) = \phi(h) $ if $h \in H$ and $f_{1}(g) =0$ if $g \not\in H$.

If, following Lemma \ref{5.7}, ${\mathcal G}$ is a set of coset representatives for 
$H \backslash G  $ then a $k$-basis for $c-\underline{{\rm Ind}}_{H}^{G}(\phi)$ is given by 
 \[   \{  g \cdot f_{1}  \ | \  g \in {\mathcal G} \}  . \]
 
 For $g \in G$ we have
\[ \begin{array}{ll}
(g \cdot f_{1})(x) = f_{1}(xg) & =  \left\{  \begin{array}{cc}
0 & {\rm if } \ xg \not\in H, \\
\\
\phi(xg)  & {\rm if } \ xg \in H, \\
\end{array} \right.   \\
\\
& =  \left\{  \begin{array}{cc}
0 & {\rm if } \  x  \not\in Hg^{-1}, \\
\\
\phi(xg)  & {\rm if } \ x \in Hg^{-1}. \\
\end{array} \right. 
\end{array} \]

Before going further let us introduce the presence of $(H, \phi)$ into the notation.
 \end{em}
\end{example}
\begin{definition}
\label{5.12}
\begin{em}

Let $H$ be a closed subgroup of $G$ containing the centre, $Z(G)$, which is compact open modulo $Z(G)$. Let $\phi : H \longrightarrow k^{*}$ be a continuous character of $H$. Denote by $X_{c}(H, \phi)$ the $k$-vector space of functions $f: G \longrightarrow k  $ such that

(i) \  $f(hg) = \phi(h)f(g)$ for all $h \in H, g \in G$,

(ii)  \  there is a compact open subgroup $K_{f} \subseteq G$ such that $f(gk) = f(g)$ for all $g \in G, k \in K_{f}$,

(ii)  $f$ is compactly supported modulo $H$.

As in \S\ref{5.6}, the left action of $G$ on  $X_{c}(H, \phi)$ is given by $(g \cdot f)(x)= f(xg)$ and therefore
\[   \Sigma :  G  \longrightarrow  {\rm Aut}_{k}( X_{c}(H, \phi))  \]
gives a smooth representation of $G$ - denoted by $\Sigma = c-{\rm Ind}_{H}^{G}(\phi)$.

Henceforth we shall denote the map written as $f_{1}$ in Example \ref{5.11} by 
\linebreak
$f_{(H, \phi)} \in 
X_{c}(H, \phi)$.

Therefore, for $g \in G$, we have
\[ \begin{array}{ll}
(g \cdot f_{(H, \phi)})(x) = f_{(H, \phi)}(xg) & =  \left\{  \begin{array}{cc}
0 & {\rm if } \ xg \not\in H, \\
\\
\phi(xg)  & {\rm if } \ xg \in H, \\
\end{array} \right.   \\
\\
& =  \left\{  \begin{array}{cc}
0 & {\rm if } \  x  \not\in Hg^{-1}, \\
\\
\phi(xg)  & {\rm if } \ x \in Hg^{-1}. \\
\end{array} \right. 
\end{array} \]
\end{em}
\end{definition}
\begin{definition}
\label{5.13}
\begin{em}

For $(H, \phi)$ and $(K, \psi)$ as in Definition \ref{5.12},  write 
\linebreak
$[(K, \psi), g, (H, \phi)]$ for any triple consisting of $g \in G$, characters $\phi, \psi$ on 
subgroups $H, K \leq G$, respectively such that 
\[    (K, \psi) \leq (g^{-1}Hg, (g)^{*}(\phi)) \]
which means that $K \leq  g^{-1}Hg$ and that $\psi(k) = \phi(h)$ where $k = g^{-1}hg$ for $h \in H, k \in K$.

Let ${\mathcal H}$ denote the $k$-vector space with basis given by these triples. Define a product on these triples by the formula 
\[  [(H, \phi), g_{1}, (J, \mu)]  \cdot  [(K, \psi), g_{2}, (H, \phi)] =   [(K, \psi), g_{1}g_{2}, (J, \mu)]  \]
and zero otherwise. This product makes sense because 

(i)  \   if $K \leq g_{2}^{-1} H g_{2}$ and
 $H \leq g_{1}^{-1} J g_{1}$ then $K \leq  g_{2}^{-1} H g_{2} \leq   g_{2}^{-1} g_{1}^{-1} J g_{1} g_{2} $ 
 
 and 
 
 (ii)  \  if $\psi(k) = \phi(h) = \mu(j), $ where $k = g_{2}^{-1}hg_{2},  h = g_{1}^{-1}j g_{1}$ then
 \linebreak
  $k = g_{2}^{-1}  g_{1}^{-1}j g_{1}   g_{2}$. 
 
 This product is clearly associative and we define an algebra ${\mathcal H}_{cmc}(G)$ to be ${\mathcal H}$ modulo the relations (c.f. Appendix \S4)
 \[   [(K, \psi), gk, (H, \phi)]  = \psi(k^{-1}) [(K, \psi), g, (H, \phi)]  \]
and
 \[     [(K, \psi), hg, (H, \phi)]  = \phi(h^{-1}) [(K, \psi), g, (H, \phi)] .    \]
 
 We observe that
 \[  [(K, \psi), g, (H, \phi)] = [( g^{-1}Hg, g^{*}\phi), g , (H, \phi)] \cdot [(K, \psi), 1, ( g^{-1}Hg, g^{*}\phi)]\]
 
We shall refer to this algebra as the compactly supported modulo the centre (CSMC-algebra) of $G$.
\end{em}
\end{definition}

\begin{lemma}{$_{}$}
\label{5.14}
\begin{em}

Let $[(K, \psi), g, (H, \phi)]$ be a triple as in Definition \ref{5.13}. Associated to this triple define a left $k[G]$-homomorphism
\[  [(K, \psi), g, (H, \phi)] : X_{c}(K, \psi) \longrightarrow X_{c}(H, \phi) \]
by the formula $g_{1} \cdot f_{(K, \psi)} \mapsto (g_{1}g^{-1}) \cdot  f_{(H, \phi) }$.
\end{em}
\end{lemma}

For a proof, which is the same as in the case when $G$ is finite, can be found in (the Appendix on induction in the case of finite groups).

\begin{theorem}{$_{}$}
\begin{em}
\label{5.15}

Let ${\mathcal M}_{c}(G)$ denote the partially order set of pairs $(H, \phi)$ as in Definitions \ref{5.12} and \ref{5.13} (so that $ X_{c}(H, \phi) = c-{\rm Ind}_{H}^{G}(\phi)$). Then, when each $n_{\alpha} =1$, 
\[        M_{c}(\underline{n}, G) = \oplus_{\alpha \in {\mathcal A}, (H, \phi) \in {\mathcal M}_{c}(G)}  \\underline{n}_{\alpha} X_{c}(H, \phi)  \]
is a left $k[G] \times {\mathcal H}_{cmc}(G)$-module. For a general distribution of multiplicities
$ \{ n_{\alpha} \}$ it is Morita equivalent to a left $k[G] \times {\mathcal H}_{cmc}(G)$-module.
\end{em}
\end{theorem}

{\bf Proof}

We have only to verify associativity of the module multiplication, which is obvious.  $\Box$

\begin{definition}{$_{k[G]}{\bf mon}$, the monomial category of $G$}
\label{5.16}
\begin{em}

The monomial category of $G$ is the additive category (non-abelian) whose objects are the $k$-vector spaces given by direct sums of $ X_{c}(H, \phi)$'s of \S\ref{5.15} and whose morphisms are elements of the hyperHecke algebra ${\mathcal H}_{cmc}(G)$. In other words the subcategory of the category of 
$k[G] \times {\mathcal H}_{cmc}(G)$-modules of which one example is $ M_{c}(\underline{n}, G)$ in \S\ref{5.15}.
\end{em}
\end{definition}

\section{Idempotented algebras  (\cite{DB96} p.309)}
\begin{definition}
\label{6.1}
\begin{em}

Let $k$ be a field and $H$ a $k$-algebra. Let ${\mathcal E}$ denote a set of idempotents of $H$. Assume that if $e_{1}, e_{2} \in {\mathcal E}$ then there exists $e_{0} \in {\mathcal E}$ such that $e_{0}e_{1} = e_{1}e_{0} = e_{1}$ and $e_{0}e_{2} = e_{2}e_{0} = e_{2}$. In addition assume for every $\phi \in H$ that there exists $e \in {\mathcal E}$ such that $e \phi = \phi e = \phi$. 

With these assumptions $H$ is called an idempotented $k$-algebra.

Write $f \leq e$ if $ef=fe=f$. This gives ${\mathcal E}$ the structure of a partially ordered set (i.e. a poset).

If $R$ is a ring and $e$ an idempotent denote $eRe$ by $R[e]$. If $M$ is a left $R$-module write $M[e]$ for the $R[e]$-module $eM$. If $H$ is an idempotented algebra then $H[e]$ is a $k$-algebra with unit $e$ and $M[e]$ is an $H[e]$-module.

$M$ is smooth if $M = \bigcup_{e \in {\mathcal E}} \ M[e]$ and is admissible if it is smooth and for each $e \in {\mathcal E}$ we have ${\rm dim}_{k}(M[e]) < \infty$.

 If $(H_{i}, {\mathcal E}_{i})$ are idempotented algebras for $i=1,2$ then so is
 $H_{1} \otimes H_{2}$ with idempotents $e_{1} \otimes e_{2}$ for $e_{i} \in {\mathcal E}_{i}$.
 \end{em}
 \end{definition}
 \begin{dummy}{The idempotented algebra ${\mathcal H}_{cmc}(G)$ }
\label{6.2}
 \begin{em}

 Let ${\mathcal E}$ be the collection of finite additive combinations in ${\mathcal H}_{cmc}(G)$, the algebra of Definition \ref{5.13}, of the form
 \[  e =  \sum_{i=1}^{n} \ [(H_{i}, \phi_{i}), 1 , (H_{i}, \phi_{i})]   \] 
 in which $(H_{i}, \phi_{i}) = (H_{j}, \phi_{j})$ if and only if $i=j$. Then $e \cdot e = e$ and all the idempotents in ${\mathcal H}_{cmc}(G)$ have this form. 
 
  We shall write $e_{(H, \phi)}$ for the idempotent $[(H, \phi), 1 , (H, \phi)]$.
 \end{em}
 \end{dummy}
 
 Define the homomorphism
\[  [(K', \psi'), g, (H', \phi')]  : X_{c}(K, \psi) \longrightarrow  X_{c} (H, \phi)   \]
to be zero unless $K', \psi') = (K, \psi)$ and $(H, \phi) =  (H', \phi')$. The following result is clear.
 \begin{theorem}{$_{}$}
\label{6.3}
 \begin{em}
  
(i) \    In \S\ref{6.2} $({\mathcal H}_{cmc}(G), {\mathcal E})$ is an idempotented algebra and $M_{c}(G)$ is an  ${\mathcal H}_{cmc}(G)$-module in the category of smooth $k[G]$-modules.
  
(ii) \    In this idempotented algebra  $e =  \sum_{i=1}^{n} \ [(H_{i}, \phi_{i}), 1 , (H_{i}, \phi_{i})] $ and $f$ satisfy $ef =fe =f$ in and only if the idempotent $f$ is a subsum of $e$, which fits very nicely with the $f \leq e$ notation. 

(iii) \  If $M_{c}(\underline{n}, G)$ is the module of Theorem \ref{5.15} then 
$M_{c}(\underline{n}, G)[e]$ is the direct sum of the $\underline{n}_{\alpha} X_{c}(H, \phi)$'s for which $e_{(H, \phi)}$ appears in the sum for $e$. 
 \end{em}
 \end{theorem} 
 \begin{dummy}{Hecke algebras}
\label{6.4}
 \begin{em}
 
 The Hecke algebra of a locally compact, totally disconnected group is a related idempotented algebra.
 \end{em}
 \end{dummy}
Let $G$ be a locally compact, totally disconnected group. Assume that $G$ is unimodular - that is, the left invariant Haar measure equals the right-invariant Haar measure of $G$ (\cite{DB96} p.137).
  
 The Hecke algebra of $G$, denoted by ${\mathcal H}_{G}$ is the space $C_{c}^{\infty}(G)$ of locally constant,compactly supported $k$-valued functions on $G$ with the convolution product (\cite{DB96} p.140 and p.255)
 \[   (\phi_{1} * \phi_{2})(g)  =  \int_{G} \phi_{1}(gh) \phi_{2}(h^{-1}) dh  = \int_{G}  \phi_{1}(h) \phi_{2}(h^{-1}g) dh  . \]
 This integral requires only one of $\phi_{1}, \phi_{2}$ to be compactly supported in order to land in ${\mathcal H}_{G}$.
 
 Suppose that $K_{0} \subseteq G$ is a compact, open subgroup. Define an idempotent
 \[   e_{K_{0}} = \frac{1}{{\rm vol}(K_{0})} \cdot \chi_{K_{0}} \]
 where $\chi_{K_{0}}$ is the characteristic function of $K_{0}$. If $K_{0} \subseteq K_{1}$ then $e_{K_{0}} * e_{K_{1}} = e_{K_{1}} $.
 
 This is seen using left invariance of the Haar measure
 \[  \int_{G} \  \frac{\chi_{K}(zh)}{{\rm vol}(K)} \frac{\chi_{H}(h^{-1})}{{\rm vol}(H)} =  
  \int_{G} \  \frac{\chi_{K}(h)}{{\rm vol}(K)} \frac{\chi_{H}(h^{-1}z)}{{\rm vol}(H)} . \]
  The integrand is zero unless $h \in K$ and then it is zero unless $z \in H$. When $z \in H$ we are integrating
\[    \int_{G} \  \frac{\chi_{K}(h)}{{\rm vol}(K)} \frac{1}{{\rm vol}(H)} =  \frac{\chi_{H}(z)}{{\rm vol}(H)},  \]
as required.

 ${\mathcal H}_{G}$ is an idempotented algebra because $G$ has a base of neighbourhoods consisting of compact open subgroups.
 
 A function $f \in  {\mathcal H}_{G}$ is called $K$-finite if the subspace spanned by all its (left) translates by $K$ is finite-dimensional (\cite{DB96} p.299).
 
   \section{Monomial morphisms as convolution products}
 
It is my belief and eventual intention that the material of this section will remain true for the general $G$ as in \S2 provided that all continuous $k$-valued characters on compact, open subgroups have finite image. 

However, throughout  this section I shall assume that $G$ is a locally profinite group whose centre $Z(G)$ is compact. Let $H$ be a subgroup which is compact, open modulo the centre. Let $k$ be an algebraically closed field for which all continuous characters $\phi : H \longrightarrow k^{*}$ have finite image when $H$ is compact, open. 

The following two results give some examples of $G$ for which $Z(G)$ is compact.
\begin{lemma}{$_{}$}
\label{7.1}
\begin{em}

Let $K$ be a $p$-adic local field. Then $Z(SL_{n}K)$ is finite. In particular it is compact.
\end{em}
\end{lemma}

{\bf Proof}

Consider the relation
\[ \begin{array}{l}
 \left( \begin{array}{ccccc}
x_{1} & 0 &  \cdot &  \cdots & 0 \\
0 & x_{2} &  \cdot &  \cdots & 0 \\
0&0& \cdots  & \cdots & 0 \\
\vdots & \vdots &   \vdots  & \vdots & \vdots \\
0&0& \cdots  & x_{n-1} & 0 \\
0&0&\cdots  & 0 & x_{n} \\
\end{array}  \right) 
 \left( \begin{array}{ccccc}
a_{1,1} & a_{1,2} &  \cdot &  \cdots & a_{1,n} \\
a_{2,1} & a_{2,2} &  \cdot &  \cdots & a_{2,n} \\
a_{3,1}&a_{3,2}& \cdots  & \cdots & a_{3,n} \\
\vdots & \vdots &   \vdots  & \vdots & \vdots \\
a_{n-1,1}&a_{n-1,2} & \cdots  & a_{n-1,n-1} & a_{n-1,n} \\
a_{n,1}&a_{n,2}&\cdots  & a_{n,n-1} & a_{n,n} \\
\end{array}  \right)   \\
\\
=  \left( \begin{array}{ccccc}
a_{1,1} & a_{1,2} &  \cdot &  \cdots & a_{1,n} \\
a_{2,1} & a_{2,2} &  \cdot &  \cdots & a_{2,n} \\
a_{3,1}&a_{3,2}& \cdots  & \cdots & a_{3,n} \\
\vdots & \vdots &   \vdots  & \vdots & \vdots \\
a_{n-1,1}&a_{n-1,2} & \cdots  & a_{n-1,n-1} & a_{n-1,n} \\
a_{n,1}&a_{n,2}&\cdots  & a_{n,n-1} & a_{n,n} \\
\end{array}  \right) 
 \left( \begin{array}{ccccc}
x_{1} & 0 &  \cdot &  \cdots & 0 \\
0 & x_{2} &  \cdot &  \cdots & 0 \\
0&0& \cdots  & \cdots & 0 \\
\vdots & \vdots &   \vdots  & \vdots & \vdots \\
0&0& \cdots  & x_{n-1} & 0 \\
0&0&\cdots  & 0 & x_{n} \\
\end{array}  \right) 
\end{array}  \]
In the $(i,j)$ entry we find $x_{i} a_{i,j} = a_{i,j}x_{j}$ and since we may suppose $a_{i,j} \not= 0$ we see that $x_{1} = x_{2} = \ldots = x_{n}$ and $x_{1}^{n} = 1$. Therefore $Z(SL_{n}K) = \mu_{n}(K)$, the group of $n$-th roots of unity in $K$. $\Box$

\begin{lemma}{$_{}$}
\label{7.2}
\begin{em}

Let $K$ be a $p$-adic local field with ring of integers ${\mathcal O}_{K}$ and prime $\pi_{K}$. Then $Z(GL_{n}K/ \langle \pi_{K} \rangle) \cong {\mathcal O}_{K}^{*}$ . In particular it is compact. Here
$\langle \pi_{K} \rangle$ denotes the centre subgroup generated by $\pi_{K}$ times the identity matrix.
\end{em}
\end{lemma}

{\bf Proof}

The relation used in the proof of \S\ref{7.1} implies that for each $(i,j)$ we have
$\pi_{K}^{\alpha} x_{i} a_{i,j} = a_{i,j}x_{j} \pi_{K}^{\beta}$ for some pair $\alpha, \beta$. Therefore we may suppose that $x_{1} \in {\mathcal O}_{K}^{*}$ and that $x_{j} = x_{1} \pi_{K}^{e_{j}}$. Now taking a matrix with $a_{1,j} a_{j,1} \not= 0$ for $j=2,3, \ldots , n$ we find that
$\pi_{K}^{\alpha} x_{1}  = x_{j} \pi_{K}^{\beta} = x_{1}  \pi_{K}^{\beta + e_{j}}$ for $j=2,3, \ldots , n$. This implies that $x_{1} \pi_{K}^{e} = x_{2} = x_{3} = \ldots = x_{n}$ which implies that $e=0$. $\Box$

The next two results ensure that we are free to use convolution products in our context.
\begin{lemma}{$_{}$}
\label{7.3}
\begin{em}

Let $G$ be a locally profinite group whose centre $Z(G)$ is compact. If $H$ is a subgroup of $G$, containing $Z(G)$, which is compact, open modulo the centre then $H$ is compact, open. 
\end{em}
\end{lemma} 

{\bf Proof}

The is a compact open subset $C$ of $G$ such that $H = Z(G) \cdot C$. Multiplication is a continuous map from the compact space $Z(G) \times C$ onto $H$ so that $H$ is compact. Furthermore any point of $H$ may be written as $h = z \cdot c$ with $z \in Z(G)$ and $c \in C$. Therefore $z \cdot N \subseteq H$
for any open neighbourhood of $c$ in $C$ is an open neighbourhood of $h$ in $H$, which is therefore open. $\Box$
\begin{lemma}{$_{}$}
\label{7.4}
\begin{em}

Let $G$ be a locally profinite group whose centre $Z(G)$ is compact and let $H$ be a subgroup which is compact, open modulo the centre. Let $k$ be an algebraically closed field for which all continuous characters $\phi : H \longrightarrow k^{*}$ have finite image when $H$ is compact, open. Then the vector space, $X_{c}$, of \S\ref{5.6} on which $c-{\rm Ind}_{H}^{G}(k_{\phi})$ is defined is a subspace of the Hecke algebra of $G$,  ${\mathcal H}_{G}$, the space of locally constant,compactly supported $k$-valued functions on $G$.
\end{em}
\end{lemma}

{\bf Proof:}

By \S\ref{5.7} it suffices to verify that the function $f_{w}$ of \S\ref{5.6} is locally constant, compactly supported for $w = 1 \in k^{*}$. This function is given by the formula
\[   f_{1}(x)  =  \left\{  \begin{array}{cc}
0 & {\rm if } \ x \not\in H, \\
\\
\phi(x)  & {\rm if } \ x \in H, \\
\end{array} \right.     \]
By \S\ref{7.3} $H$, the support of $f_{1}$, is compact. Since the image of $\phi$ is finite the function $f_{1}$ is locally constant. $\Box$

  Recall from \S\S\ref{5.13}-\ref{5.14} that we have defined
\[  [(K, \psi), g, (H, \phi)] : X_{c}(K, \psi) \longrightarrow X_{c}(H, \phi) \]
by the formula $g_{1} \cdot f_{(K, \psi)} \mapsto (g_{1}g^{-1}) \cdot  f_{(H, \phi) }$.

If $\chi_{W}$ is the characteristic function of $W \subseteq G$ we may define $g_{1} \cdot f_{(K, \psi)}$ using characteristic functions in the following manner. By definition
\[  g_{1} \cdot f_{(K, \psi)}(x) = \left\{ \begin{array}{ll}
\psi(xg_{1}^{-1}) & {\rm if} \ xg_{1}^{-1} \in K,  \\
\
0 & {\rm if} \ xg_{1}^{-1}  \not\in K  .
\end{array} \right.   \]
Suppose that $v_{1}, \ldots , v_{t}$ are coset representatives for $K/{\rm Ker}(\psi)$. Then, if $xg_{1}^{-1} \in K$ we must have $xg_{1}^{-1} \in {\rm Ker}(\psi)v_{j(xg_{1}^{-1} )}$ for some $1 \leq j(xg_{1}^{-1} ) \leq t$ and therefore 
$\psi(xg_{1}^{-1} ) = \psi(v_{j(xg_{1}^{-1} )})$. Hence we have the fomula
\[  g_{1} \cdot f_{(K, \psi)} = \sum_{j=1}^{t} \ \psi(v_{j}) \chi_{{\rm Ker}(\psi) v_{j} g_{1}}  \]
because $\bigcup {\rm Ker}(\psi) v_{j} g_{1} = Kg_{1}$ so that the right hand side is zero unless 
\linebreak
$xg_{1}^{-1} \in K$ and is $\psi(v_{j_{0}})$ precisely when $j_{0} = j(xg_{1}^{-1})$.

Next, from Definition \ref{5.13}
\[    (K, \psi) \leq (g^{-1}Hg, (g)^{*}(\phi)) \] 
implies that $\psi(k) = \phi(h)$ where $k = g^{-1}hg$ for $h \in H, k \in K$. Therefore if $k \in {\rm Ker}(\psi)$ then $h \in {\rm Ker}(\phi)$ and so ${\rm Ker}(\psi) \leq g^{-1} {\rm Ker}(\phi)g$.

Consider the convolution product
\[   \chi_{g_{1}{\rm Ker}(\psi)}  * \chi_{g^{-1}{\rm Ker}(\phi)}(z) =  \int_{G}   \chi_{g_{1}{\rm Ker}(\psi)}(h)
 \chi_{g^{-1}{\rm Ker}(\phi)}(h^{-1}z)  dh .   \]
 The integrand is zero unless $h \in g_{1}{\rm Ker}(\psi)$ in addition to the condition 
 \[  z \in  hg^{-1}{\rm Ker}(\phi)  \subseteq  g_{1}{\rm Ker}(\psi)g^{-1}{\rm Ker}(\phi) = g_{1}g^{-1}  g{\rm Ker}(\psi)g^{-1}{\rm Ker}(\phi) \subseteq  g_{1}g^{-1}  {\rm Ker}(\phi) \]
 and conversely. Therefore
 \[   \chi_{g_{1}{\rm Ker}(\psi)}  * \chi_{g^{-1}{\rm Ker}(\phi)} = {\rm vol}( g_{1}{\rm Ker}(\psi)) \chi_{ g_{1}g^{-1}  {\rm Ker}(\phi)} . \]
 
 Similarly, if $v \in K$ and $u \in H$, we have a convolution product
 \[  \chi_{g_{1}{\rm Ker}(\psi)v}  * \chi_{g^{-1}{\rm Ker}(\phi)u}(z) =  \int_{G}   \chi_{g_{1}{\rm Ker}(\psi)v}(h)
 \chi_{g^{-1}{\rm Ker}(\phi)u}(h^{-1}z)  dh .     \]
 The integrand is zero unless $h \in g_{1}{\rm Ker}(\psi)v$ in addition to the condition
 \[ z \in hg^{-1}{\rm Ker}(\phi)u  \subseteq  g_{1}{\rm Ker}(\psi)v g^{-1}{\rm Ker}(\phi)u
 \subseteq   g_{1}g^{-1}  {\rm Ker}(\phi) (gvg^{-1}) \cdot u \]
 and conversely. Therefore
 \[   \chi_{g_{1}{\rm Ker}(\psi)v}  * \chi_{g^{-1}{\rm Ker}(\phi)u} = {\rm vol}(g_{1}{\rm Ker}(\psi)v) 
 \chi_{ g_{1}g^{-1} {\rm Ker}(\phi) gvg^{-1}u }  . \]
 
 \begin{lemma}{$_{}$}
\label{7.5}
 \begin{em}
 
 Suppose that $v_{1}, \ldots , v_{t} \in K$ is a set of coset representatives for $K/{\rm Ker}(\psi)$.
 Then
 \[ g_{1} \cdot f_{(K, \psi)} = \sum_{j=1}^{t} \ \psi(v_{j}) \cdot  \chi_{{\rm Ker}(\psi)v_{j}g_{1}^{-1} } .  \]
 \end{em}
 \end{lemma}
 
 {\bf Proof:}
 
 Consider the functions in the equation applied to $x \in G$. The left hand side is zero if $xg_{1} \not\in K$ which is equivalent to there being no $j$ such that $xg_{1} \in {\rm Ker}(\psi)v_{j}$ or
 $x \in {\rm Ker}(\psi)v_{j}g_{1}^{-1}$. Under these conditions every characteristic function on the right hand side also vanishes on $x$. On the other hand if $xg_{1} \in K$ there exists a unique $j_{0}$ such that $x \in  {\rm Ker}(\psi)v_{j_{0}}g_{1}^{-1} $ and so, evaluated at $xg_{1}$, there is one and only one term on the right hand side which contributes. It yields $\psi(v_{j_{0}})$ which is the value of $g_{1} \cdot f_{(K, \psi)}$ at $x$, as required. $\Box$
 
 \begin{dummy}
\label{7.6}
 \begin{em}
 
 The image $\phi(H)$ is a finite cyclic group, being a finite subgroup of $k^{*}$, which contains $\phi(gKg^{-1}) = \psi(K)$. Therefore there exist $v_{1}, \ldots , v_{t}$ which are coset representatives for $K/{\rm Ker}(\psi)$ and $u_{1}, \ldots , u_{s}$ which give distinct cosets in $H/{\rm Ker}(\phi)$ such that the set $\{  (gv_{i}g^{-1}) u_{j} \ | \ 1 \leq i \leq t, \ 1 \leq j \leq s \}$ is a set of coset representatives for $H/{\rm Ker}(\phi)$.
 \end{em}
 \end{dummy}
 \begin{definition}
\label{7.7}
 \begin{em}
  Define an involution $T :  C_{c}^{\infty}(G) \longrightarrow C_{c}^{\infty}(G)$ by 
  \linebreak
  $T(F)(x) = F(x^{-1})$. For example $T( \chi_{{\rm Ker}(\psi)v_{j}g_{1}^{-1} }) =   \chi_{ g_{1}{\rm Ker}(\psi)v_{j}^{-1} }$.
  
  In the notation of \S\ref{7.6} set
  \[ \Phi_{[(K, \psi), g, (H, \phi)]} = \sum_{j=1}^{s} \phi(u_{j}) \cdot \chi_{g^{-1} {\rm Ker}(\phi) u_{j} } . \]
 \end{em}
 \end{definition}
 
 \begin{theorem}{$_{}$}
\label{7.8}
 \begin{em}
 
 In the notation of Definition \ref{7.7}
 \[   \begin{array}{ll}
   [(K, \psi), g, (H, \phi)]( g_{1} \cdot  f_{(K, \psi)}) & = g_{1}g^{-1} \cdot  f_{(H, \phi)} \\
   \\
   & =   \frac{1}{ {\rm vol}(Ker(\psi)) }  T(   T(g_{1} \cdot f_{(K, \psi)})  *   \Phi_{[(K, \psi), g, (H, \phi)]} ) .
 \end{array}  \]
 \end{em}
 \end{theorem}
 
 {\bf Proof:}  
 
 We observe that  $\psi(v_{i})( \chi_{{\rm Ker}(\psi)v_{i}g_{1}^{-1} })(x) = \psi(v_{i}) = \psi( xg_{1})$
 if $x \in {\rm Ker}(\psi)v_{i}g_{1}^{-1} =  v_{i}{\rm Ker}(\psi)g_{1}^{-1}$ and zero otherwise. Therefore 
 \[ T(\psi(v_{i})( \chi_{{\rm Ker}(\psi)v_{i}g_{1}^{-1} }))(x) = \psi(v_{i})( \chi_{{\rm Ker}(\psi)v_{i}g_{1}^{-1} })(x^{-1}) =  \psi(v_{i}) \]
 if $x^{-1} \in {\rm Ker}(\psi)v_{i}g_{1}^{-1} $ and zero otherwise. In the non-zero case $x \in g_{1} {\rm Ker}(\psi) v_{i}^{-1}$ and $\psi(v_{i}) = \psi(g_{1}^{-1}x)^{-1}$ so that
\[   T(\psi(v_{i})( \chi_{{\rm Ker}(\psi)v_{i}g_{1}^{-1} }))  = 
 \psi(v_{i})^{-1} \chi_{g_{1} {\rm Ker}(\psi) v_{i}^{-1}}     .  \]
 
  From Lemma \ref{7.1} we have 
 \[  \begin{array}{l}
  T(g_{1} \cdot f_{(K, \psi)})  
  =   \sum_{i=1}^{t} \ \psi(v_{i})^{-1} \cdot   \chi_{ g_{1} {\rm Ker}(\psi)v_{i}^{-1} } .
  \end{array} \]
Therefore
\[  \begin{array}{l}
   T(g_{1} \cdot f_{(K, \psi)})  *   \Phi_{[(K, \psi), g, (H, \phi)]} \\
   \\
   =   \sum_{i=1}^{t} \  \sum_{j=1}^{s} \  \psi(v_{i})^{-1}  \phi(u_{j})  (  \chi_{ g_{1} {\rm Ker}(\psi)v_{i}^{-1} } *   \chi_{g^{-1} {\rm Ker}(\phi) u_{j} })   \\
   \\
   =  \sum_{i=1}^{t} \  \sum_{j=1}^{s} \  \psi(v_{i})^{-1}  \phi(u_{j}) {\rm vol}({\rm Ker}(\psi)) 
   \chi_{ g_{1}g^{-1} {\rm Ker}(\phi) (gv_{i}^{-1}g^{-1})u_{j}}   .
 \end{array} \]
 Hence
 \[ \begin{array}{l}
  T(   T(g_{1} \cdot f_{(K, \psi)})  *   \Phi_{[(K, \psi), g, (H, \phi)]} )  \\
  \\
  = T( \sum_{i=1}^{t} \  \sum_{j=1}^{s} \  \psi(v_{i})^{-1}  \phi(u_{j}) {\rm vol}({\rm Ker}(\psi)) 
   \chi_{ g_{1}g^{-1} {\rm Ker}(\phi) (gv_{i}^{-1}g^{-1})u_{j}} ) \\
   \\
    = T( \sum_{i=1}^{t} \  \sum_{j=1}^{s} \  \phi(gv_{i}g^{-1})^{-1}  \phi(u_{j}) {\rm vol}({\rm Ker}(\psi)) 
   \chi_{ g_{1}g^{-1} {\rm Ker}(\phi) (gv_{i}^{-1}g^{-1})u_{j}} ) \\
   \\
       =  {\rm vol}({\rm Ker}(\psi))   \sum_{i=1}^{t} \  \sum_{j=1}^{s} \ T( \phi((gv_{i}^{-1}g^{-1}) u_{j}) 
   \chi_{ g_{1}g^{-1} {\rm Ker}(\phi) (gv_{i}^{-1}g^{-1})u_{j}} ) \\
   \\
      =  {\rm vol}({\rm Ker}(\psi))   \sum_{i=1}^{t} \  \sum_{j=1}^{s} \  \phi((gv_{i}^{-1}g^{-1}) u_{j})^{-1} 
   \chi_{    u_{j}^{-1}  (gv_{i}g^{-1})   {\rm Ker}(\phi) gg_{1}^{-1} }  \\
   \\
   = {\rm vol}(Ker(\psi))  g_{1}g^{-1} \cdot  f_{(H, \phi)} ,
 \end{array} \]
 by Lemma \ref{7.1}. $\Box$
 \begin{remark}
 \label{7.9}
 \begin{em}
 
 (i) \  Theorem \ref{7.8} has shown that, under the special conditions which were stated at the start of this section, the morphisms of the monomial category $_{k[G]}{\bf mon}$ of Definition \ref{5.16} are given in terms of the convolution product of \S\ref{6.4} of the Hecke algebra ${\mathcal H}_{G}$.
 
(ii) \  My belief is that Theorem \ref{7.8} remains true in general, in some sense, providing that all continuous characters $\phi : H \longrightarrow k^{*}$ have finite image when $H$ is compact, open.  
This belief is based on the following: \cite{Sn18} claims to construct for each admissible representation $V$ of $G$ a monomial resolution in the derived category $_{k[G]}{\bf mon}$\footnote{In a later section I shall give a self-contained construction of these resolutions based on the hyperHecke algebra and which applies to any  $V$ is  ${\mathcal M}_{cmc, \underline{\phi}}(G)$-admissible $V$.} and (see \S9; also \cite{PD84} pp.2-3) such $V$ are intimately related to Hecke modules. Therefore one should expect a connection between the morphisms in that resolution and convolutions products.

The difficulty, in the case of a general locally profinite group $G$, with the treatment of this section is that 
$X_{c}(H, \phi)$'s are spaces of locally constant functions which are compactly supported modulo $H$, rather than actually being compactly supported. 

It might be that I can get away with using the Schwartz space of locally constant, compactly supported functions of $G/Z(G)$, but I have not yet had time to examine this generalisation\footnote{To that end, as a novice, I should re-read \S9, several sections of \cite{DB96} on Hecke modules and the 
material of (\cite{AJS79} \S1.11 p.63).}.
 \end{em}
 \end{remark}
 
\section{The bar-monomial resolution: I. finite modulo the centre case}
\begin{dummy}
\label{8.1.1}
\begin{em}

Let $G$ is a finitely generated group with centre $Z(G)$ and finite quotient group $G/Z(G)$ with other notation as in \S2 and other earlier sections. In this section I summarise (\cite{Sn18} Chapter One, \S2).

Denote by $\hat{G}$ the group of character homomorphisms ${\rm Hom}(G, k^{*})$ from $G$ to $k^{*}$. For $H$ be a subgroup of $G$ containing $Z(G)$ denote by $\hat{H}_{\underline{\phi}}$ the finite subset of $\hat{H}$ of characters which are equal to $\underline{\phi} $ when restricted to $Z(G)$.

For an arbitrary $k$-algebra $A$ we write $_{A}{\bf mod}$ (resp. ${\bf mod}_{A}$) for the category of left (resp. right) $A$-modules. We denote by $_{A}{\bf lat}$ (resp. ${\bf lat}_{A}$) the category of left (resp. right) $A$-lattices i.e. the subcategory of $_{A}{\bf mod}$ consisting of those $A$-modules which are finitely generated and $A$-projective. The rank of a free $k$-module $M$ will be denoted by 
${\rm rk}_{k}(M)$. When $A = k[G]$ we have subcategories  $_{k[G], \underline{\phi}}{\bf mod} \subset   \   _{k[G]}{\bf mod}$ and $_{k[G], \underline{\phi}}{\bf lat} \subset   \   _{k[G]}{\bf lat}$ whose objects are those on which $Z(G)$ acts via $\underline{\phi}$.

Recall that the poset of ${\mathcal M}_{\underline{\phi}}(G)$ of pairs $(H, \phi)$ admits a left $G$-action by conjugation for which the $G$-orbit of $(H, \phi)$ will be denoted by $(H, \phi)^{G}$.
\end{em}
\end{dummy}

\begin{definition}
\label{8.1.2}
\begin{em}

 A finite $(G, \underline{\phi})$-lineable left $k[G]$-module $M$\footnote{Here I have taken my own terminological advice given in the footnote to (\cite{Sn18} Chapter One, Definition 1.2). } is a left $k[G]$-module together with a fixed finite direct sum decomposition
 \[ M = M_{1} \oplus \cdots \oplus M_{m} \]
 where each of the $M_{i}$ is a free $k$-module of rank one on which $Z(G)$ acts via $\underline{\phi}$
 and the $G$-action permutes the $M_{i}$. The $M_{i}$'s are called the lines of $M$. For $1 \leq i \leq m$ let $H_{i}$ denote the subgroup of $G$ with stabilises the line $M_{i}$. Then there exists a unique $\phi_{i} \in \hat{H_{i}}_{\underline{\phi}}$ such that $h \cdot v = \phi_{i}(h) v$ for all $v \in M_{i}, h \in H_{i}$. The pair $(H_{i}, \phi_{i}) \in {\mathcal M}_{ \underline{\phi}}(G)$ is called the stabilising pair of $M_{i}$.
 
 The $k$-submodule of $M$ given by
 \[   M^{((H, \phi))} = \oplus_{1 \leq i \leq m, \ (H, \phi) \leq (H_{i}, \phi_{i})} \ M_{i}  \]
 is called the $(H, \phi)$-fixed points of $M$.
 
 A morphism between $(G, \underline{\phi})$-lineable modules from $M$ to $N = N_{1} \oplus \cdots \oplus N_{n}$ is defined to be a $k[G]$-module homomorphism $f: M \longrightarrow N$ such that
 \[ f(  M^{((H, \phi))} ) \subseteq N^{((H, \phi))}   \] 
 for all $(H, \phi) \in {\mathcal M}_{\underline{\phi}}(G)$. 
 
 The (left) finite $(G, \underline{\phi})$-lineable modules and their morphisms define an additive category denoted by $_{k[G], \underline{\phi}}{\bf mon}$.
 
 By definition each $(G, \underline{\phi})$-lineable module is a $k$-free $k[G]$-module so there is a forgetful functor 
 \[  {\mathcal V} :  \  _{k[G], \underline{\phi}}{\bf mon} \longrightarrow  \     _{k[G], \underline{\phi}}{\bf mod}  .  \]
 \end{em}
\end{definition} 
\begin{remark}{Some natural operations and constructions}
\label{8.1.3}
\begin{em}

There are several operations which are obvious lifts to $_{k[G], \underline{\phi}}{\bf mon}$ of well-known operations in $_{k[G], \underline{\phi}}{\bf mod}$. This means that the resulting functors commute with the forgetful functor $ {\mathcal V} :  \  _{k[G], \underline{\phi}}{\bf mon} \longrightarrow   \     _{k[G], \underline{\phi}}{\bf mod}  $.

Direct sum and tensor product, denoted by $M \oplus N$ and $M \otimes N$, are lineable in the obvious way as is the $k[G]$-module ${\rm Hom}_{k}(M,N)$, giving an object of  $_{k[G], \underline{\phi_{M}} \cdot \underline{\phi_{N}}^{-1}  }{\bf mon}$. A homomorphism, preserving centres, $f: G' \longrightarrow  G$ induces a restriction map ${\rm Res}_{f}$  from  $_{k[G], \underline{\phi}}{\bf mon}  $ to $ _{k[G'], \underline{\phi} \cdot f}{\bf mon}  $. 

 Suppose that $Z(G) \subseteq H \subseteq G$ and the index of $H$ in $G$ is finite and 
 $P \in   _{k[H], \underline{\phi}}{\bf mon}  $ then the usual induced $k[G]$-module construction (see, Appendix: Comparison of Inductions) ${\rm Ind}_{H}^{G}(P)$ gives a lineable module 
$\underline{{\rm Ind}}_{H}^{G}(P)$ in such a way that, if the stabilising pair for $P_{i}$ is 
$(H_{i}, \phi_{i})$, then the stabilising pair of the Line $g \otimes_{k[H]} P_{i}$ is $g(H_{i}, \phi_{i})$. 
 
Analogues of the usual distributivity isomorphism of direct sums over tensor products, the Frobenius reciprocity isomorphism and the Mackey decomposition isomorphism 
all hold in $ _{k[G], \underline{\phi}}{\bf mon}$.

We call a $M$ in $_{k[G], \underline{\phi}}{\bf mon}$ indecomposable if it is not isomorphic to a non-trivial direct sum $N \oplus P$ in $_{k[G], \underline{\phi}}{\bf mon}$. 
\end{em}
\end{remark}
\begin{proposition}{$_{}$}
\label{8.1.4}
\begin{em}
The set of $(G, \underline{\phi})$-lineeable modules given by 
 \[ \{ \underline{{\rm Ind}}_{H}^{G}(k_{\phi}) \ | \ (H, \phi) \in G \backslash {\mathcal M}_{ \underline{\phi}}(G) \} \]
 is a full set of pairwise non-isomorphic representatives for the isomorphism classes of indecomposable objects in $_{k[G], \underline{\phi}}{\bf mon}$. Moreover any object in $_{k[G], \underline{\phi}}{\bf mon}$
is canonically isomorphic to the direct sum of objects in this set.
\end{em}
\end{proposition}
\begin{dummy}
\label{8.1.5}
\begin{em}

Let $[(K, \psi) , g , (H, \phi)]$ be one of the basic generators of the hyperHecke algebra 
${\mathcal H}_{cmc})(G)$ of \S2 then we have a morphism
\[  [(K, \psi) , g , (H, \phi)]  \in   {\rm Hom}_{  _{k[G], \underline{\phi}}{\bf mon}}( \underline{{\rm Ind}}_{K}^{G}(k_{\psi}) , \underline{{\rm Ind}}_{H}^{G}(k_{\phi}))  \]
defined by the same formula as in the case of induced modules (see, Appendix: Comparison of Inductions). In addition the composition of morphisms in $_{k[G], \underline{\phi}}{\bf mon}$ coincides with the product in the hyperHecke algebra.

 Let $(K, \psi)   \in  {\mathcal M}_{\underline{\phi}}(G)$ and let $N$ be an object of  $_{k[G], \underline{\phi}}{\bf mon}$. Then there is a $k$-linear isomorphism
\[  {\rm Hom}_{  _{k[G], \underline{\phi}}{\bf mon}}( \underline{{\rm Ind}}_{K}^{G}(k_{\psi}) , N) \stackrel{\cong}{\longrightarrow}  N^{((K, \psi))}  \]
given by $f \mapsto f(1 \otimes_{K} 1)$. The inverse isomorphism is given by 
\[ n \mapsto ((g \otimes_{K} v \mapsto v g \cdot n))   . \]
\end{em}
\end{dummy}
\begin{lemma}{Projectivity in $_{k[G], \underline{\phi}}{\bf mon}$}
\label{8.1.6}
\begin{em}

Consider the diagram
\[   M  \stackrel{h}{\longrightarrow}  N  \stackrel{f}{\longleftarrow}  P  \]
in which $M, P \in  _{k[G], \underline{\phi}}{\bf mon}$ and $N \in _{k[G], \underline{\phi}}{\bf mod}$ with $h,f$ being morphisms in $_{k[G], \underline{\phi}}{\bf mod}$.  Assume, for all $(H, \phi) \in {\mathcal M}_{\underline{\phi}}(G)$, that
\[   f(P^{((H, \phi))}) \subseteq h(M^{((H, \phi))})    . \]
Then there exists $j \in  {\rm Hom}_{  _{k[G], \underline{\phi}}{\bf mon}}( P , M)$ such that $h \cdot j = f$.

In particular we include the situation where $N' \in _{k[G], \underline{\phi}}{\bf mon}$ with $h, f$ being morphisms to $N'$ in $_{k[G], \underline{\phi}}{\bf mon}$ and the diagram above being the result of applying the forgetful functor ${\mathcal V}$ with $N = {\mathcal V}(N')$.
\end{em}
\end{lemma}
\begin{dummy}
\label{8.1.7}
\begin{em}
Let $V$ be a finitely generated $k[G]$-module on which $Z(G)$ acts via the central character $\underline{\phi}$. That is, $V$ is an object of 
$_{k[G], \underline{\phi}}{\bf mod}$. I shall define the notion of a 
$_{k[G], \underline{\phi}}{\bf mon}$-resolution of $V$. This is a chain complex of morphisms in 
$_{k[G], \underline{\phi}}{\bf mon}$ with certain properties which will ensure that it is exists and is unique up to chain homotopy in $_{k[G], \underline{\phi}}{\bf mon}$.

For $V \in _{k[G], \underline{\phi}}{\bf mod}$ and $(H, \phi) \in  {\mathcal M}_{\underline{\phi}}(G)$ define the $(H, \phi)$-fixed points of $V$ by
\[   V^{(H, \phi)} = \{ v \in V \ | \ h \cdot v = \phi(h)v \ {\rm for \ all} \ h \in H\}.\]
Clearly $g(V^{(H, \phi)}) = V^{g(H, \phi)}$, $V^{(Z(G), \underline{\phi})} = V$ and $(K, \psi) \leq (H, \phi)$ implies that $ V^{(H, \phi)}  \subseteq  V^{(K, \psi)} $. Note that $f \in {\rm Hom}_{_{k[G], \underline{\phi}}{\bf mod}}(V,W)$ satisfies $f( V^{(H, \phi)}) \subseteq  W^{(H, \phi)}$ for all $(H, \phi) \in  {\mathcal M}_{\underline{\phi}}(G)$. In addition, if $M \in  _{k[G], \underline{\phi}}{\bf mon}$ then $M^{((H, \phi))} \subseteq M^{(H, \phi)}$ so that 
\end{em}
\end{dummy}
\begin{definition}{(\cite{Sn18} Chapter One \S2)}
\label{8.1.8}
\begin{em}

Let $V \in _{k[G], \underline{\phi}}{\bf mod}$. A $_{k[G], \underline{\phi}}{\bf mon}$-resolution of $V$ is a chain complex
\[   M_{*}:  \hspace{20pt}  \ldots \stackrel{\partial_{i+1}}{\longrightarrow}  \ M_{i+1} \stackrel{\partial_{i}}{\longrightarrow} M_{i} \stackrel{\partial_{i-1}}{\longrightarrow}  \ldots   \stackrel{\partial_{1}}{\longrightarrow}  \ M_{1} \stackrel{\partial_{0}}{\longrightarrow} M_{0}   \]
with $M_{i} \in _{k[G], \underline{\phi}}{\bf mon}$ and $\partial_{i} \in {\rm Hom}_{_{k[G], \underline{\phi}}{\bf mon}}(M_{i+1}, M_{i})$ for all $i \geq 0$ together with $\epsilon \in  {\rm Hom}_{_{k[G], \underline{\phi}}{\bf mod}}({\mathcal V}(M_{0}), V)$ such that
\[  \ldots  \stackrel{\partial_{i}}{\longrightarrow} M_{i}^{((H, \phi))} \stackrel{\partial_{i-1}}{\longrightarrow}  \ldots   \stackrel{\partial_{1}}{\longrightarrow}  \ M_{1}^{((H, \phi))} \stackrel{\partial_{0}}{\longrightarrow} M_{0}^{((H, \phi))}   \stackrel{\epsilon}{\longrightarrow}  V^{(H, \phi)} \longrightarrow  0 \]
is an exact sequence of $k$-modules for each $(H, \phi) \in {\mathcal M}_{\underline{\phi}}(G)$. In particular, when $(H, \phi) = (Z(G), \underline{\phi})$ we see that
\[  \ldots  \stackrel{\partial_{i}}{\longrightarrow} M_{i} \stackrel{\partial_{i-1}}{\longrightarrow}  \ldots   \stackrel{\partial_{1}}{\longrightarrow}  \ M_{1} \stackrel{\partial_{0}}{\longrightarrow} M_{0}   \stackrel{\epsilon}{\longrightarrow}  V \longrightarrow  0 \]
is an exact sequence in $_{k[G], \underline{\phi}}{\bf mod}$.
\end{em}
\end{definition}
\begin{proposition}{$_{}$}
\label{8.1.9}
\begin{em}

  Let $V \in _{k[G], \underline{\phi}}{\bf mod}$ and let
\[ \ldots \longrightarrow  M_{n} \stackrel{\partial_{n-1}}{\longrightarrow}  M_{n-1} \stackrel{\partial_{n-2}}{\longrightarrow} \ldots  \stackrel{\partial_{0}}{\longrightarrow}  M_{0}  \stackrel{\epsilon}{\longrightarrow} V  \longrightarrow 0 \]
be a $_{k[G], \underline{\phi}}{\bf mon}$-resolution of $V$. Suppose that
\[ \ldots \longrightarrow  C_{n} \stackrel{\partial'_{n-1}}{\longrightarrow}  C_{n-1} \stackrel{\partial'_{n-2}}{\longrightarrow} \ldots  \stackrel{\partial'_{0}}{\longrightarrow}  C_{0}  \stackrel{\epsilon'}{\longrightarrow} V  \longrightarrow 0 \]
a chain complex where each $\partial'_{i}$ and $C_{i}$ belong to $_{k[G], \underline{\phi}}{\bf mon}$
 and $\epsilon'$ is a $_{k[G], \underline{\phi}}{\bf mod}$ homomorphism such that $\epsilon'(C_{0}^{((H, \phi))}) \subseteq V^{(H, \phi)}$ for each $(H, \phi) \in {\mathcal M}_{\underline{\phi}}(G)$. 
 
 Then there exists a chain map of $_{k[G], \underline{\phi}}{\bf mon}$-morphisms  $\{ f_{i} : C_{i} \longrightarrow  M_{i} , i \geq 0 \}$ such that
\[   \epsilon \cdot f_{0} = \epsilon', \     f_{i-1} \cdot \partial'_{i}  = \partial_{i} \cdot f_{i} \  {\rm for \ all} \ i \geq 1 .\]

In addition, if $\{ f'_{i} : C_{i} \longrightarrow  M_{i} , i \geq 0 \}$ is another chain map of $_{k[G], \underline{\phi}}{\bf mon}$-morphisms such that $ \epsilon \cdot f_{0} =  \epsilon \cdot f'_{0}$ then there exists a $_{k[G], \underline{\phi}}{\bf mon}$-chain homotopy
 $\{  s_{i} : C_{i} \longrightarrow M_{i+1} ,  \  {\rm for \ all} \ i \geq 0  \}$
 such that $\partial_{i} \cdot s_{i} + s_{i-1} \cdot \partial'_{i} = f_{i} - f'_{i}$ for all $i \geq 1$ and $f_{0}- f'_{0} = \partial_{0} \cdot s_{0}$.
\end{em}
\end{proposition}
\vspace{2pt}

{\bf Proof}
\vspace{2pt}

This is the usual homological algebra argument using Lemma \ref{8.1.6}. $ \Box$  
\begin{remark}
\label{8.1.10}
\begin{em}

(i) \ Needless to say, Proposition \ref{8.1.9} has an analogue to the effect that every $_{k[G], \underline{\phi}}{\bf mod}$-homomorphism $V \longrightarrow V'$ extends to a 
$_{k[G], \underline{\phi}}{\bf mon}$-morphism between the monomial resolutions of $V$ and $V'$, if they exist, and the extension is unique up to $_{k[G], \underline{\phi}}{\bf mon}$-chain homotopy.

(ii) \  The category $_{k[G], \underline{\phi}}{\bf mon}$ is additive but not abelian. Homological algebra (e.g. a projective resolution) is more conveniently accomplished in an abelian category. To overcome this difficulty we shall embed $_{k[G], \underline{\phi}}{\bf mon}$ into more convenient abelian categories. This is reminiscent of the Freyd-Mitchell Theorem which embeds every abelian category into a category of modules.
\end{em}
\end{remark}
\begin{dummy}{The functor category $funct_{k}^{o}(_{k[G], \underline{\phi}}{\bf mon}, _{k}{\bf mod})$}
\label{8.1.11}
\begin{em}

Let $funct_{k}^{o}(_{k[G], \underline{\phi}}{\bf mon}, _{k}{\bf mod})$ denote the category of contravariant  functors, ${\mathcal F}, {\mathcal G}$ etc, 
from $_{k[G], \underline{\phi}}{\bf mon}$ to the category of finitely generated $k$-modules whose morphisms are 
$k$-linear natural transformations $\alpha : {\mathcal F} \longrightarrow  {\mathcal G}$ etc.

Let $_{k[G], \underline{\phi}}{\bf mod}$ denote the category of finite rank $k[G]$-modules with central character $\underline{\phi}$ (see \S\ref{8.1.1}). Consider the functor
\[ {\mathcal I} :   _{k[G], \underline{\phi}}{\bf mod}  \longrightarrow  funct_{k}^{o}(_{k[G], \underline{\phi}}{\bf mon}, _{k}{\bf mod})   \]
given on objects by
\[ {\mathcal I}(V) =  {\rm Hom}_{_{k[G], \underline{\phi}}{\bf mod} }({\mathcal V}(-), V)  \]
with ${\mathcal I}(\beta : V \longrightarrow W) = (f \mapsto \beta \cdot f)$.

The relation of $ {\mathcal I}$ to morphisms, particularly to the morphisms given by the 
$[(K, \psi) , g , (H, \phi)]$'s of the hyperHecke algebra, is analysed in detail in \cite{Sn18} Chapter One \S3).
\end{em}
\end{dummy}
\begin{proposition}{The functor ${\mathcal I}$ on morphisms}
\label{8.1.12}
\begin{em}

Let $V,W  \in _{k[G], \underline{\phi}}{\bf mod}$  and set ${\mathcal F}(-) = {\mathcal I}(V)(-)$, ${\mathcal G}(-) = {\mathcal I}(W)(-)$. Given a natural transformation $\alpha : {\mathcal F} \longrightarrow  {\mathcal G}$ for each $M  \in _{k[G], \underline{\phi}}{\bf mon}$ we have
\[ \alpha(M) :  {\mathcal F}(M) \longrightarrow  {\mathcal G}(M) \] such that if $\beta : M \longrightarrow M'$ is a morphism in $_{k[G], \underline{\phi}}{\bf mon}$ we have a commutative diagram

If we have a homomorphism $\gamma : V \longrightarrow W$ in $_{k[G], \underline{\phi}}{\bf mod}$ we obtain a natural transformation 
\[ \gamma_{*} : {\mathcal F}  \longrightarrow  {\mathcal G} \]
 given by $\gamma_{*}(M)(f) = \gamma \cdot f  \in {\mathcal G}(M)$ for all $f: M \longrightarrow V$ in ${\mathcal F}(M)$.
 
 Any natural transformation 
\[ \alpha : {\mathcal F}  \longrightarrow  {\mathcal G} \]
 there is a unique homomorphism $\gamma : V \longrightarrow W$ in $_{k[G], \underline{\phi}}{\bf mod}$ such that $\alpha = \gamma_{*}$, which completely determines $\alpha$.
\end{em}
\end{proposition}
\begin{proposition}{$_{}$}
\label{8.1.13}
\begin{em}
 Let ${\mathcal I}$ denote the functor of \S\ref{8.1.11} and define a functor
 \[  {\mathcal J} :   _{k[G], \underline{\phi}}{\bf mon}  \longrightarrow  funct_{k}^{o}(_{k[G], \underline{\phi}}{\bf mon}, _{k}{\bf mod}) \]
by ${\mathcal J}(M) =  {\rm Hom}_{ _{k[G], \underline{\phi}}{\bf mon} }(-, M) $. 

Then the category 
$funct_{k}^{o}(_{k[G], \underline{\phi}}{\bf mon}, _{k}{\bf mod})$ is abelian. Furthermore both ${\mathcal I}$ and ${\mathcal J}$ are full embeddings (i.e. bijective on morphisms and hence injective on isomorphism classes of objects).
\end{em}
\end{proposition}
\vspace{2pt}
\begin{proposition}
\label{8.1.14}
\begin{em}

For $M \in _{k[G], \underline{\phi}}{\bf mon}$ the functor ${\mathcal J}(M)$ in 
\linebreak
$funct_{k}^{o}(_{k[G], \underline{\phi}}{\bf mon}, _{k}{\bf mod}) $ is projective.
\end{em}
\end{proposition}
\begin{definition}
\label{8.1.15}
\begin{em}
 
 Let $M \in _{k[G], \underline{\phi}}{\bf mon}, V \in _{k[G], \underline{\phi}}{\bf mod}$. Define a $k$-linear isomorphism $ {\mathcal K}_{M,V} $ of the form
 \[   {\rm Hom}_{_{k[G], \underline{\phi}}{\bf mod}}({\mathcal V}(M), V) \stackrel{ {\mathcal K}_{M,V} }{\longrightarrow}
 {\rm Hom}_{funct_{k}^{o}(_{k[G], \underline{\phi}}{\bf mon}, _{k}{\bf mod})}({\mathcal J}(M), {\mathcal I}(V)) \]
by sending $f : {\mathcal V}(M)  \longrightarrow  V$ to the natural transformation 
\[  {\mathcal K}_{M,V}(N)  :  {\mathcal J}(M)(N) \longrightarrow    {\mathcal I}(V)(N)  \]
given by $h \mapsto f \cdot {\mathcal V}(h)$  for all $N \in _{k[G], \underline{\phi}}{\bf mon}$
\[  {\rm Hom}_{ _{k[G], \underline{\phi}}{\bf mon} }(N, M)  \longrightarrow 
 {\rm Hom}_{_{k[G], \underline{\phi}}{\bf mod}}({\mathcal V}(N), V) .   \]
 
 The inverse isomorphism is given by $ {\mathcal K}_{M,V}^{-1}(\phi) = \phi(M)(1_{M}) $ where $1_{M}$ denotes the identity morphism on $M$. 
 
 In fact ${\mathcal K}$ is a functorial equivalence of the form
  \[  {\mathcal K} :  {\rm Hom}_{_{k[G], \underline{\phi}}{\bf mod}}({\mathcal V}(-), -) \stackrel{ \cong }{\longrightarrow}
 {\rm Hom}_{funct_{k}^{o}(_{k[G], \underline{\phi}}{\bf mon}, _{k}{\bf mod})}({\mathcal J}(-), {\mathcal I}(-)) \]
\end{em}
\end{definition}
\begin{theorem}{$_{}$}
\label{8.1.16}
\begin{em}

Let
\[  \ldots  \stackrel{\partial_{i}}{\longrightarrow} M_{i} \stackrel{\partial_{i-1}}{\longrightarrow}  \ldots   \stackrel{\partial_{1}}{\longrightarrow}  \ M_{1} \stackrel{\partial_{0}}{\longrightarrow} M_{0}   \stackrel{\epsilon}{\longrightarrow}  V \longrightarrow  0 \]
be a chain complex with $M_{i} \in _{k[G], \underline{\phi}}{\bf mon}$ for $i \geq 0$, $V \in _{k[G], \underline{\phi}}{\bf mod}$,
\linebreak
 $\partial_{i} \in {\rm Hom}_{ _{k[G], \underline{\phi}}{\bf mon} }(M_{i+1}, M_{i})$  and $\epsilon \in  {\rm Hom}_{_{k[G], \underline{\phi}}{\bf mod}}({\mathcal V}(M_{0}), V)$. Then the following are equivalent:
 
(i)  \ $M_{*} \longrightarrow V$ is a $_{k[G], \underline{\phi}}{\bf mon}$-resolution of $V$.

(ii) \ The sequence
\[ \begin{array}{l}
 \ldots  \stackrel{{\mathcal J}(\partial_{i})}{\longrightarrow} {\mathcal J}(M_{i}) \stackrel{{\mathcal J}(\partial_{i-1})}{\longrightarrow}  \ldots   \stackrel{{\mathcal J}(\partial_{1})}{\longrightarrow}  \ {\mathcal J}(M_{1}) \stackrel{{\mathcal J}(\partial_{0})}{\longrightarrow} {\mathcal J}(M_{0})  
  \stackrel{  {\mathcal K}_{M_{0},V}(\epsilon)  }{\longrightarrow}  {\mathcal I}(V) \longrightarrow  0
\end{array} \]
is exact in $funct_{k}^{o}(_{k[G], \underline{\phi}}{\bf mon}, _{k}{\bf mod})$.
\end{em}
\end{theorem}
\begin{remark}
\label{8.1.17}
\begin{em}

Theorem \ref{8.1.16} together with Proposition \ref{8.1.13} and Proposition \ref{8.1.14} imply that the map
\[  ( M_{*} \stackrel{\epsilon}{\longrightarrow} V )   \ \mapsto  \
 ( {\mathcal J}(M_{*}) \stackrel{{\mathcal K}_{M_{0},V}(\epsilon)}{\longrightarrow}  {\mathcal I}(V)) \]
 is a bijection between the $_{k[G], \underline{\phi}}{\bf mon}$-resolutions of $V$ and the
 projective resolutions of ${\mathcal I}(V)$ consisting of objects from the subcategory
 ${\mathcal J}(_{k[G], \underline{\phi}}{\bf mon})$. 
\end{em}
\end{remark}

We are now going to pass from functors to a category of modules.
\begin{dummy}{The functor $\Phi_{M}$}
\label{8.1.18}
\begin{em}

Let $M \in  _{k[G], \underline{\phi}}{\bf mon}$ and let ${\mathcal A}_{M} = {\rm Hom}_{_{k[G], \underline{\phi}}{\bf mon}}(M,M)$, the ring of endomorphisms on $M$ under composition. In the present context 
${\mathcal A}_{M}$ is a finitely generated $k$-algebra.
\end{em}
\end{dummy}

I shall show that there is an equivalence of categories between 
\linebreak
$funct_{k}^{o}(_{k[G], \underline{\phi}}{\bf mon}, _{k}{\bf mod})$ and the category of right modules ${\bf mod}_{{\mathcal A}_{M}}$ for a suitable choice of $M$.

We have a functor
\[ \Phi_{M} : funct_{k}^{o}(_{k[G], \underline{\phi}}{\bf mon}, _{k}{\bf mod}) \longrightarrow  {\bf mod}_{{\mathcal A}_{M}} \] 
given by $\Phi({\mathcal F}) = {\mathcal F}(M)$. Right multiplication by $ z \in {\mathcal A}_{M}$ on $v \in  {\mathcal  F}(M)$ is given by 
\[ v \# z =  {\mathcal F}(z)(v) \]
where $ {\mathcal F}(z) : {\mathcal F}(M) \longrightarrow  {\mathcal F}(M)$ is the left $k$-module morphism obtained by applying ${\mathcal F}$ to the endomorphism $z$. This is a right-${\mathcal A}_{M}$ action since
 \[ v \# (zz_{1})  =  {\mathcal F}(zz_{1})(v) = ({\mathcal F}(z_{1}) \cdot   {\mathcal F}(z))(v) = 
 {\mathcal F}(z_{1})( {\mathcal F}(z)(v)) = (v  \# z)   \# z_{1} . \]
 
 In the other direction define a functor
 \[ \Psi_{M} : {\bf mod}_{{\mathcal A}_{M}}  \longrightarrow  funct_{k}^{o}(_{k[G], \underline{\phi}}{\bf mon}, _{k}{\bf mod}) ,  \]
for $P \in {\bf mod}_{ {\mathcal A}_{M}}$, by
\[ \Psi_{M}(P) = {\rm Hom}_{{\mathcal A}_{M}}( {\rm Hom}_{_{k[G], \underline{\phi}}{\bf mon}}(M, - ) , P) . \]
Here, for $N \in  _{k[G], \underline{\phi}}{\bf mon}$, ${\rm Hom}_{_{k[G], \underline{\phi}}{\bf mon}}(M, - )$ is a right ${\mathcal A}_{M}$-module via pre-composition by endomorphisms of $M$. For a homomorphism of ${\mathcal A}_{M}$-modules 
\linebreak
$f : P  \longrightarrow Q$ the map $\Psi_{M}(f)$ is given by composition with $f$.

Next we consider the composite functor 
\[   \Phi_{M} \cdot \Psi_{M} :   {\bf mod}_{{\mathcal A}_{M}}   \longrightarrow 
 {\bf mod}_{{\mathcal A}_{M}}.  \]
This is given by $P \mapsto {\rm Hom}_{{\mathcal A}_{M}}( {\rm Hom}_{_{k[G], \underline{\phi}}{\bf mon}}(M, M) , P) = {\rm Hom}_{{\mathcal A}_{M}}( {\mathcal A}_{M} , P) $ so that there is an obvious natural transformation $\eta : 1 \stackrel{\cong}{\longrightarrow}  \Phi_{M} \cdot \Psi_{M}$ such that $\eta(P)$ is an isomorphism for each module $P$.

Now consider the composite functor
\[   \Psi_{M} \cdot \Phi_{M} :   funct_{k}^{o}(_{k[G], \underline{\phi}}{\bf mon}, _{k}{\bf mod})  \longrightarrow  funct_{k}^{o}(_{k[G], \underline{\phi}}{\bf mon}, _{k}{\bf mod}) .  \]
For a functor ${\mathcal F}$ we shall define a natural transformation 
\[ \epsilon_{{\mathcal F}} : {\mathcal F} \longrightarrow {\rm Hom}_{{\mathcal A}_{M}}( {\rm Hom}_{_{k[G], \underline{\phi}}{\bf mon}}(M, - ) ,  {\mathcal F}(M)) =    \Psi_{M} \cdot \Phi_{M}({\mathcal F}) . \]
For $N \in  _{k[G], \underline{\phi}}{\bf mon}$ we define
\[ \epsilon_{{\mathcal F}}(N)   : {\mathcal F}(N)   \longrightarrow {\rm Hom}_{{\mathcal A}_{M}}( {\rm Hom}_{_{k[G], \underline{\phi}}{\bf mon}}(M, N ) ,  {\mathcal F}(M))   \]
by the formula $v \mapsto (f \mapsto {\mathcal F}(f)(v))$.

\begin{theorem}{$_{}$}
\label{8.1.19}
\begin{em}

Let $S \in  _{k[G], \underline{\phi}}{\bf mon}$ be the finite $(G, \underline{\phi})$-lineable $k$-module given by
\[ S = \oplus_{(H, \phi) \in  {\mathcal M}_{\underline{\phi}}(G)} \ \underline{{\rm Ind}}_{H}^{G}(k_{\phi})  . \]
Then, in the notation of \S\ref{8.1.18}, 
\[ \Phi_{S} : funct_{k}^{o}(_{k[G], \underline{\phi}}{\bf mon}, _{k}{\bf mod}) \longrightarrow  {\bf mod}_{{\mathcal A}_{S}} \] 
and 
 \[ \Psi_{S} : {\bf mod}_{{\mathcal A}_{S}}  \longrightarrow  funct_{k}^{o}(_{k[G], \underline{\phi}}{\bf mon}, _{k}{\bf mod})   \]
 are inverse equivalences of categories. In fact, the natural transformations $\eta$ and $\epsilon$ are isomorphisms of functors when $M= S$. 
\end{em}
\end{theorem}
\begin{remark}
\label{8.1.20}
\begin{em}

Theorem \ref{8.1.19} is true when $S$ is replaced by any $M$ which is the direct sum of $\underline{{\rm Ind}}_{H}^{G}(k_{\phi}) $'s containing at least one pair $(H, \phi)$ from each $G$-orbit of $ {\mathcal M}_{\underline{\phi}}(G)$. That is, for any $(G, \underline{\phi})$-lineable $k$-module containing 
\[  \oplus_{(H, \phi) \in  G \backslash {\mathcal M}_{\underline{\phi}}(G)} \ \underline{{\rm Ind}}_{H}^{G}(k_{\phi})  \]
as a summand. This remark is established by Morita theory.
\end{em}
\end{remark}

\begin{dummy}{The bar resolution}
\label{8.1.21}
\begin{em}

We begin this section by recalling the two-sided bar-resolution for $A$-modules. Let $k$ be a commutative Noetherian ring and let $A$ be a (not necessarily commutative) $k$-algebra. For each integer $p \geq 0$ set
\[ B_{p}(M, A, N) =  M \otimes_{k} A \otimes_{k} A \otimes_{k} \ldots \otimes_{k} A \otimes_{k} N \]
in which there are $p$ copies of $A$ and $M \in   {\bf mod}_{A},  N \in   _{A}{\bf mod}$.
Define 
\[ d : B_{1}(M, A, N) \longrightarrow  B_{0}(M, A, N)  \]
by $d(m \otimes a \otimes n) = m \cdot a \otimes n  -   m \otimes a \cdot n $. For $p \geq 2$ define
\[ d : B_{p}(M, A, N) \longrightarrow  B_{p-1}(M, A, N)  \]
by
\[  \begin{array}{l}
d(m \otimes a_{1} \otimes \ldots \ \otimes a_{p} \otimes n) \\
\\
= 
m \cdot a_{1} \otimes \ldots \ \otimes a_{p} \otimes n  \\
\\
+ \sum_{i=1}^{p-1} (-1)^{i} m \otimes a_{1} \otimes \ldots \otimes  a_{i} \cdot a_{i+1} \otimes \ldots \otimes a_{p} \otimes n  \\
\\
+ (-1)^{p}  m \otimes a_{1} \otimes \ldots \ \otimes a_{p} \cdot n .
\end{array}  \]
Setting $N=A$ we define
\[ \epsilon : B_{0}(M, A, A) = M \otimes_{k} A \longrightarrow  M \]
by $\epsilon(m \otimes a) = m \cdot a$.
\end{em}
\end{dummy}
 The following result is well-known.
{\begin{proposition}{$_{}$}
\label{8.1.22}
\begin{em}

In the situation of \S\ref{8.1.21} the chain complex
\[ \begin{array}{l}
   \ldots  \longrightarrow  B_{p}(M, A, A) \stackrel{d}{\longrightarrow}   B_{p-1}(M, A, A) \stackrel{d}{\longrightarrow} 
\ldots  \\
\\
\hspace{20pt}  \ldots \longrightarrow B_{1}(M, A, A) \stackrel{d}{\longrightarrow}  B_{0}(M, A, A) \stackrel{\epsilon}{\longrightarrow} M 
\longrightarrow  0  
\end{array} \]
is a free right-$A$-module resolution of $M$.
\end{em}
\end{proposition} 
\begin{dummy}
\label{8.1.23}
\begin{em}

As in \S\ref{8.1.15} and \S\ref{8.1.18}, let $M \in  _{k[G], \underline{\phi}}{\bf mon}, V \in _{k[G], \underline{\phi}}{\bf mod}$ and let ${\mathcal A}_{M} = {\rm Hom}_{_{k[G], \underline{\phi}}{\bf mon}}(M,M)$, the ring of endomorphisms on $M$ under composition.
For $i \geq 0$ define $\tilde{M}_{M,i} \in  \   _{k}{\bf mod}$ by ($i$ copies of $ {\mathcal A}_{M}$)
\[ \tilde{M}_{M,i} = {\rm Hom}_{_{k[G], \underline{\phi}}{\bf mod}}({\mathcal V}(M), V) \otimes_{k}  {\mathcal A}_{M} \otimes_{k} \ldots \otimes_{k}  {\mathcal A}_{M}  \]
and set
\[ \underline{M}_{M,i} =  \tilde{M}_{M,i} \otimes_{k}   {\rm Hom}_{_{k[G], \underline{\phi}}{\bf mon}}(-, M)  . \]
Hence $ \underline{M}_{M,i}  \in  funct_{k}^{o}(_{k[G], \underline{\phi}}{\bf mon}, _{k}{\bf mod}) $ and in fact the values of this functor are not merely objects in $_{k}{\bf mod}$ because they have a natural right $ {\mathcal A}_{M}$-module structure, defined as in \S\ref{8.1.18}.

If $i \geq 1$ we defined natural transformations $d_{M,0}, d_{M,1}, \ldots , d_{M,i}$ in the following way.
Define
\[ d_{M,0} : \underline{M}_{M,i}   \longrightarrow  \underline{M}_{M,i-1} \] 
by
\[ d_{M,0}( f \otimes \alpha_{1} \otimes \ldots \otimes \alpha_{i} \otimes u) = f(- \cdot \alpha_{1}) \otimes \alpha_{2} \ldots \otimes \alpha_{i} \otimes u .\]
The map $f(- \cdot \alpha_{1}) : {\mathcal V}(M) \longrightarrow V$ is a 
$_{k[G], \underline{\phi}}{\bf mod}$- homomorphism since $\alpha_{i}$ acts on the right of $M$. 

For $1 \leq j \leq i-1$ we define
\[ d_{M,j}  :  \underline{M}_{M,i}  \longrightarrow   \underline{M}_{M,i-1} \] 
by
\[ d_{M,j}( f \otimes \alpha_{1} \otimes \ldots \otimes \alpha_{i} \otimes u) =
f \otimes \alpha_{1} \ldots \otimes \alpha_{j} \alpha_{j+1} \otimes \ldots \otimes \alpha_{i} \otimes u . \]
Finally 
\[ d_{M,i}  :   \underline{M}_{M,i}  \longrightarrow   \underline{M}_{M,i-1} \] 
is given by
\[  d_{i}(M)( f \otimes \alpha_{1} \otimes \ldots \otimes \alpha_{i} \otimes u) =  
f \otimes \alpha_{1} \otimes \ldots  \otimes  \alpha_{i-1} \otimes \alpha_{i} \cdot u  .\]
Since $u$ is a $_{k[G], \underline{\phi}}{\bf mon}$-morphism so is $\alpha_{i} \cdot u$ because
\[  (\alpha_{i} \cdot u)(\alpha m) = \alpha_{i}(u(\alpha m)) =  \alpha_{i}(\alpha u(m)) =  
\alpha   \alpha_{i}(u(m)) = \alpha   (\alpha_{i} \cdot u)(m)  \]
since $\alpha_{i}$ is a $_{k[G], \underline{\phi}}{\bf mon}$ endomorphism of $M$.

Next we define a natural transformation
\[ \epsilon_{M} :  \underline{M}_{M,0} \longrightarrow  {\mathcal I}(V) =  {\rm Hom}_{ _{k[G], \underline{\phi}}{\bf mod}}({\mathcal V}(-), V)  \]
by sending $f \otimes u \in  \underline{M}_{M,0}$ to $f \cdot {\mathcal V}(u) \in   {\mathcal I}(V)$.

Finally we define 
\[   d_{M} = \sum_{j=0}^{i}  \ (-1)^{j} d_{M,j} :  \underline{M}_{M,i}  \longrightarrow   \underline{M}_{M,i-1} .  \]
\end{em}
\end{dummy}
\begin{theorem}{$_{}$}
\label{8.1.24}
\begin{em}

The sequence 
\[ \ldots \stackrel{d_{M} }{\longrightarrow}  \underline{M}_{M,i}(M)  \stackrel{d_{M} }{\longrightarrow} \underline{M}_{M,i-1}(M) \ldots 
\stackrel{d_{M} }{\longrightarrow}  \underline{M}_{M,0}(M) \stackrel{\epsilon_{M} }{\longrightarrow}  {\mathcal I}(V)(M) \longrightarrow 0 \]
is the right ${\mathcal A}_{M}$-module bar resolution of  ${\mathcal I}(V)(M) $.
\end{em}
\end{theorem}
\begin{dummy}{The functorial monomial resolution of $V$}
\label{8.1.25}
\begin{em}

Let $V$ be a finite rank $k$-lattice with a left $G$-action. Let $M \in  \ _{k[G], \underline{\phi}}{\bf mon}$ and $W \in  \ _{k}{\bf lat}$. 
 Define another object $W \otimes_{k} M \in \  _{k[G], \underline{\phi}}{\bf mon}$
 by letting $G$ act only on the $M$-factor, $g(w \otimes  m) = w \otimes gm$, and defining the Lines of $W \otimes_{k} M$  to consist of the one-dimensional subspaces $\langle w \otimes L \rangle$ where $w\in W$, runs through a $k$-basis of $W$, and $L$ is a Line of $M$. Therefore, if $M'  \in  \ 
 _{k[G], \underline{\phi}}{\bf mon}$ we have an isomorphism
\[   {\rm Hom}_{ _{k[G], \underline{\phi}}{\bf mon}}(W \otimes_{k} M, M') \stackrel{\cong}{\longrightarrow} W \otimes_{k} 
{\rm Hom}_{ _{k[G], \underline{\phi}}{\bf mon} }(M, M')  \]
providing that $W$ is finite-dimensional.
Similarly we have an isomorphism
\[  {\rm Hom}_{ _{k[G], \underline{\phi}}{\bf mon}}(M, W \otimes_{k}  M') \stackrel{\cong}{\longrightarrow} W \otimes_{k} 
{\rm Hom}_{ _{k[G], \underline{\phi}}{\bf mon} }(M, M')  \]
when $W$ is finite dimensional. 

As in Theorem \ref{8.1.19}, let $S \in  _{k[G], \underline{\phi}}{\bf mon}$ be the finite $(G, \underline{\phi})$-Line Bundle over $k$ given by
\[ S = \oplus_{(H, \phi) \in  {\mathcal M}_{\underline{\phi}}(G)} \ \underline{{\rm Ind}}_{H}^{G}(k_{\phi})  . \]

As in \S\ref{8.1.23}, for $i \geq 0$ we have $\tilde{M}_{S,i} \in  \   _{k}{\bf mod}$ by ($i$ copies of $ {\mathcal A}_{S}$)
\[ \tilde{M}_{S,i} = {\rm Hom}_{_{k[G], \underline{\phi}}{\bf mod}}({\mathcal V}(S), V) \otimes_{k}  {\mathcal A}_{S} \otimes_{k} \ldots \otimes_{k}  {\mathcal A}_{S} ,  \]
which is a finite dimensional $k$-lattice. As a $k$-basis for $\tilde{M}_{S,i}$ we take the tensor product of the direct sum of bases for each $V^{(H, \phi)}$ and a basis for each $ {\mathcal A}_{S}$-factor
given by the basis of the hyperHecke algebra. 
Therefore we may form $\tilde{M}_{S,i} \otimes_{k} S \in \
_{k[G], \underline{\phi}}{\bf mon}$.

Recall from \S\ref{8.1.11} that
\[  {\rm Hom}_{_{k[G], \underline{\phi}}{\bf mod}}({\mathcal V}(S), V) = {\mathcal I}(V)(S) \cong
\oplus_{(H, \phi) \in  {\mathcal M}_{\underline{\phi}}(G)} \  V^{(H,  \phi)} . \]

We have morphisms in $_{k[G], \underline{\phi}}{\bf mon}$ for $i \geq 1$
\[  d_{0}, d_{1}, \ldots , d_{i} : \tilde{M}_{S,i} \otimes_{k} S \longrightarrow  \tilde{M}_{S,i-1} \otimes_{k} S  \]
defined on 
\[ f \otimes \alpha_{1} \otimes \ldots \otimes \alpha_{i} \otimes s \in  
{\rm Hom}_{k[G, \underline{\phi}]-mod}( {\mathcal V}(S) , V) \otimes_{k}  {\mathcal A}_{S}^{\otimes^{i}}  \otimes_{k} S \]
by
\[ \begin{array}{l}
d_{0}(f \otimes \alpha_{1} \otimes \ldots \otimes \alpha_{i} \otimes s) = 
f \cdot {\mathcal V}(\alpha_{1}) \otimes \alpha_{2} \otimes \ldots \otimes \alpha_{i} \otimes s ,   \\
\\
{\rm and \ for } \ 1 \leq j \leq i-1 \\
\\
d_{j}(f \otimes \alpha_{1} \otimes \ldots \otimes \alpha_{i} \otimes s) = 
f \otimes \alpha_{1} \otimes \ldots \otimes  \alpha_{j} \alpha_{j+1} \ldots \otimes  \alpha_{i} \otimes s  \\
\\
d_{i}(f \otimes \alpha_{1} \otimes \ldots \otimes \alpha_{i} \otimes s) = 
f \otimes \alpha_{1} \otimes \ldots \otimes \alpha_{i-1} \otimes \alpha_{i}(s) .
\end{array} \]
Setting $d = \sum_{j=0}^{i} \ (-1)^{j} d_{j}$ gives a morphism in $_{k[G], \underline{\phi}}{\bf mon}$ 
\[ d: \tilde{M}_{S,i} \otimes_{k} S \longrightarrow  \tilde{M}_{S,i-1} \otimes_{k} S \]
for all $i \geq 1$. In addition we define a homomorphism in  $_{k[G], \underline{\phi}}{\bf mod}$
\[ \epsilon :  \tilde{M}_{S,0} \otimes_{k} S = {\rm Hom}_{_{k[G], \underline{\phi}}{\bf mod}}( {\mathcal V}(S) , V) \otimes_{k} S  \longrightarrow    V \]
by $\epsilon(f \otimes s) = f(s)$.

The chain complex
\[ \begin{array}{l}
\ldots \stackrel{d}{\longrightarrow}  \tilde{M}_{S,i} \otimes_{k} S \stackrel{d}{\longrightarrow} \ldots
\stackrel{d}{\longrightarrow}    \tilde{M}_{S,1} \otimes_{k} S \stackrel{d}{\longrightarrow} 
  \tilde{M}_{S,0} \otimes_{k} S \stackrel{\epsilon}{\longrightarrow} V
\longrightarrow  0 
\end{array} \]
satisfies the conditions of Theorem \ref{8.1.16}. Therefore, abbreviating $\tilde{M}_{S,i} \otimes_{k} S$ to $M_{i}$ pro tem, this chain complex is a $_{k[G], \underline{\phi}}{\bf mon}$-resolution of $V$ if and only if the sequence
\[ \begin{array}{l}
 \ldots  \stackrel{{\mathcal J}(d)}{\longrightarrow} {\mathcal J}(M_{i}) \stackrel{{\mathcal J}(d)}{\longrightarrow}  \ldots   \stackrel{{\mathcal J}(d)}{\longrightarrow}  \ {\mathcal J}(M_{1}) \stackrel{{\mathcal J}(d)}{\longrightarrow} {\mathcal J}(M_{0})  
  \stackrel{  {\mathcal K}_{M_{0},V}(\epsilon)  }{\longrightarrow}  {\mathcal I}(V) \longrightarrow  0
\end{array} \]
is exact in $funct_{k}^{o}(_{k[G], \underline{\phi}}{\bf mon}, _{k}{\bf mod})$. By Theorem \ref{8.1.19} this chain complex of functors is exact if and only if the result of applying $\Phi_{S}$ to it is exact in the category ${\bf mod}_{{\mathcal A}_{S}}$. However, by Theorem \ref{8.1.24} with $M = S$ the resulting chain complex in ${\bf mod}_{{\mathcal A}_{S}}$ is the bar resolution, which is exact. 

We have proved the following result.
\end{em}
\end{dummy}
\begin{theorem}{Existence of the bar-monomial resolution}
\label{8.1.26}
\begin{em}

Let $k$ be a field. Then, in the notation of \S\ref{8.1.25}, The chain complex
\[ \begin{array}{l}
\ldots \stackrel{d}{\longrightarrow}  \tilde{M}_{S,i} \otimes_{k} S \stackrel{d}{\longrightarrow} \ldots
\stackrel{d}{\longrightarrow}    \tilde{M}_{S,1} \otimes_{k} S \stackrel{d}{\longrightarrow} 
  \tilde{M}_{S,0} \otimes_{k} S \stackrel{\epsilon}{\longrightarrow} V
\longrightarrow  0 
\end{array} \] 
is a $_{k[G], \underline{\phi}}{\bf mon}$-resolution of $V$.
\end{em}
\end{theorem}

The bar-monomial resolution of Theorem \ref{8.1.26} possesses a number of the usual naturality properties.

\section{The bar-monomial resolution: II. The compact, open modulo the centre case}

 Let $G$ be a locally profinite group and let $k$ be an algebraically closed field. Let $V$ be a $k$-representation of $G$ with central character $\underline{\phi}$ and that $V$ is a
 ${\mathcal M}_{cmc, \underline{\phi}}(G)$-admissible representation as in Proposition \ref{4.1}.
 
 Let  ${\mathcal H}_{cmc}(G)$ be the hyperHecke algebra, introduced in \S2. Let 
 \[ M_{c}(\underline{n}, G) = \oplus_{\alpha \in {\mathcal A}, (H, \phi) \in {\mathcal M}_{c}(G)}  \underline{n}_{\alpha} X_{c}(H, \phi)  \]
be the left $k[G] \times {\mathcal H}_{cmc}(G)$-module of Theorem \ref{5.15} form some family of strictly positive integers, $\{  \underline{n}_{\alpha} \}$.
\begin{theorem}{$_{}$}
\label{9.1.1}
\begin{em}
Replacing $S$ of Theorem \ref{8.1.19} by $M_{c}(\underline{n}, G)$ and replacing the ring ${\mathcal A}_{M}$ (when $M = S$) by ${\mathcal H}_{cmc}(G)$ we may imitate the construction of \S\ref{8.1.25}
to construct  a $_{k[G], \underline{\phi}}{\bf mon}$-resolution of $V$
\[ \begin{array}{l}
\ldots \stackrel{d}{\longrightarrow}  \tilde{M}_{M_{c}(\underline{n}, G),i} \otimes_{k} M_{c}(\underline{n}, G) \stackrel{d}{\longrightarrow} \ldots
\stackrel{d}{\longrightarrow}    \tilde{M}_{M_{c}(\underline{n}, G),1} \otimes_{k} M_{c}(\underline{n}, G) \\
\\
\hspace{50pt} \stackrel{d}{\longrightarrow} 
  \tilde{M}_{M_{c}(\underline{n}, G),0} \otimes_{k} M_{c}(\underline{n}, G) \stackrel{\epsilon}{\longrightarrow} V
\longrightarrow  0 
\end{array} \] 
\end{em}
\end{theorem} 

This result is proved using the analogues of Theorem \ref{8.1.19} and Theorem \ref{8.1.24}.
\begin{remark}
\label{9.1.2}
\begin{em}

In \cite{Sn18} this result was proved\footnote{I believe!} by reduction to the finite modulo the centre case. Also an explicit bare hands homological construction was given in the case of $GL_{2}$ of a local field.
I think that the use of the hyperHecke algebra simplifies the construction both in the compact, open modulo the centre case of this section and the general case of the next.
\end{em}
\end{remark}

\section{The monomial resolution in the general case}

Once again let $G$ be a locally profinite group and let $k$ be an algebraically closed field. Let $V$ be a $k$-representation of $G$ with central character $\underline{\phi}$ and that $V$ is a
 ${\mathcal M}_{cmc, \underline{\phi}}(G)$-admissible representation as in Proposition \ref{4.1}.
 
 First I shall recall the properties of Tammo tom Dieck's space $\underline{E}(G, {\mathcal C})$
(\cite{Sn18} Appendix IV) which is defined for a group $G$ and a family of subgroups ${\mathcal C}$
which is closed under conjugation and passage to subgroups. This space is a simplicial complex on which $G$ acts simplicially in such a way that for any subgroup $H \in {\mathcal C}$ the fixed-point set 
 $\underline{E}(G, {\mathcal C})^{H}$ is non-empty and contractible. In our case ${\mathcal C}$ will be the family of compact, open modulo the centre subgroups. 
 
 $\underline{E}(G, {\mathcal C})$ is unique up to $G$-equivariant homotopy equivalence. In the case of $GL_{n}$ of a local field, for example, the Bruhat-Tits building gives a finite-dimensional model for the tom Dieck space. 
 
 If the set of conjugacy classes maximal compact, open modulo the centre subgroups of $G$ is finite, as in the case of $GL_{n}K$ for example, one can find a local system which assigns to each compact, open  
 modulo the centre $J$ a $_{k[J], \underline{\phi}}{\bf mon}$-resolution of ${\rm Res}_{J}^{G}V$
 \[ \begin{array}{l}
\ldots \stackrel{d}{\longrightarrow}  \tilde{M}_{M_{c}(\underline{n}, J),i} \otimes_{k} M_{c}(\underline{n}, J) \stackrel{d}{\longrightarrow} \ldots
\stackrel{d}{\longrightarrow}    \tilde{M}_{M_{c}(\underline{n}, J),1} \otimes_{k} M_{c}(\underline{n}, J) \\
\\
\hspace{50pt} \stackrel{d}{\longrightarrow} 
  \tilde{M}_{M_{c}(\underline{n}, J),0} \otimes_{k} M_{c}(\underline{n}, J) \stackrel{\epsilon}{\longrightarrow} {\rm Res}_{J}^{G}V
\longrightarrow  0 .
\end{array} \] 
 
 Next one forms the double complex (\cite{Sn18} Chapter Four Theorem 3.2) given by the simplicial chain complex of the tom Dieck space in one grading and the compact, open modulo the centre 
 $_{k[J], \underline{\phi}}{\bf mon}$-resolutions in the other grading. The contribution of the resolutions corresponding to the orbit of one $J$-fixed simplex gives the compactly supported induction of that resolution.
\begin{theorem}{(\cite{Sn18} Chapter Four Theorem 3.2)}
\label{10.1.1}
\begin{em}

Let $V$ be a ${\mathcal M}_{cmc, \underline{\phi}}(G)$-admissible representation as in 
Proposition \ref{4.1}. Then the total complex of the above double complex is $_{k[G], \underline{\phi}}{\bf mon}$-resolution of $V$.
\end{em}
\end{theorem}

\section{The Bernstein centre}

Let ${\mathcal A}$ be an abelian category then its centre $Z({\mathcal A})$ is the ring of endomorphisms of the identity functor of $A$. Explicitly, for each object $A$ of $ {\mathcal A}$ there is given an endomorphism $z_{A} \in {\rm Hom}_{{\mathcal A}}(A,A)$ such that for any $f \in {\rm Hom}_{{\mathcal A}}(A,B)$ one has $z_{B}f = f z_{A}$.

If the category ${\mathcal A}$ is the product of abelian categories $( {\mathcal A}_{i})_{i \in {\mathcal I}}$
then one has  $Z({\mathcal A}) = \prod_{i \in {\mathcal I}} \ Z({\mathcal A}_{i})$.

Suppose the category ${\mathcal A}$ admits direct sums indexed by ${\mathcal I}$ such that any morphism $f: X \longrightarrow  \oplus_{i \in {\mathcal I}} \ Y_{i}$ is zero if and only if all the projections 
\[  X \stackrel{f}{\longrightarrow } \oplus_{i \in {\mathcal I}} \ Y_{i}  \stackrel{pr_{i}}{\longrightarrow }  Y_{i} \]
are zero.

This property holds for the category of algebraic (i.e. smooth) representations of a reductive group over a non-Archimedean local field (\cite{PD84} p.5).

Under the above condition ${\mathcal A}$ is the product of full subcategories ${\mathcal A}_{i}$ for $i \in {\mathcal I}$ such that

(i) \ if $X \in A_{i}$ and $Y  \in A_{j}$ then ${\rm Hom}_{{\mathcal A}}(X, Y) = 0$ if $i \not= j$ and

(ii) \ for all objects $X$  we have $X = \oplus_{i \in {\mathcal I}} \ X_{i}$ with $X_{i}$ in ${\mathcal A}_{i}$.

 \begin{definition}
 \label{11.1.1}
 \begin{em}
 
Once more let $G$ be a locally profinite group, let $k$ be an algebraically closed field and let ${\mathcal H}_{G}$ be the associated Hecke algebra as in \S\ref{6.4}. Let ${\mathcal A}$ be the category of smooth representations. Then the centre of ${\mathcal A}$ is called the Bernstein centre.

For the general linear group over a local field the Bernstein centre was determined in \cite{BZ76} and \cite{BZ77}. The determination in the general case is described in \cite{PD84}. Deligne's essay uses the description of smooth representations in term of modules over ${\mathcal H}_{G}$ (see Appendix: Smooth representations and Hecke modules) and is stated in (\cite{PD84} \S2). 
 \end{em}
 \end{definition}
 \begin{remark}
 \label{11.1.2}
\begin{em}

In this essay I cannot recall all the details of \cite{PD84}. However in this section I intend to explain an approach to the Bernstein centre of the category of smooth representation which is based on the derived category of monomial representations - that is, the derived category of modules over the hyperHecke algebra. The approach uses the monocentric conditions which were introduced in \S3.
\end{em}
\end{remark}
\begin{dummy}{Resolutions and the centre of ${\mathcal A}$}
\label{11.1.3}
\begin{em}

Suppose that ${\mathcal A}$ is an abelian category and that ${\mathcal B}$ is an additive category together with a forgetful functor $\nu : {\mathcal B} \longrightarrow {\mathcal A}$ and suppose that for each object $V \in {\rm Ob}({\mathcal A})$ we have a ${\mathcal B}$-resolution of $V$. This means a chain complex in  ${\mathcal B}$
\end{em}
\end{dummy}
\[   \stackrel{d}{\longrightarrow}  M_{i} \stackrel{d}{\longrightarrow} M_{i-1} \stackrel{d}{\longrightarrow}  \ldots \  \stackrel{d}{\longrightarrow} M_{0} \longrightarrow 0  \]
such that 
\[   \longrightarrow  \nu(M_{i}) \longrightarrow \nu(M_{i-1}) \longrightarrow  \ldots \ \longrightarrow \nu(M_{0})  \longrightarrow V \longrightarrow 0  \]
is exact in ${\mathcal A}$. In addition suppose that the association $V \mapsto M_{*}$ is functorial into the derived category of ${\mathcal B}$. Thus any two choices of ${\mathcal B}$-resolution for $V$ are chain homotopy equivalent in ${\mathcal B}$ and any morphism $ f: V \longrightarrow V'$ in ${\mathcal A}$ induces a ${\mathcal B}$-chain map, $f_{*}$ unique up to chain homotopy, between the resolutions.

Now consider a family giving an element in the centre of ${\mathcal A}$ which yields $z_{V} : V \longrightarrow V$ and $z_{V'} : V' \longrightarrow V'$ satisfying $f z_{V} = z_{V'}f$ for all $f$. Fix resolutions for $V$ and $V'$. Then $z_{V}$ induces a chain map $(z_{V})_{*}$ on $M_{*}$ and another
$(z_{V'})_{*}$ on $M'_{*}$. The morphism $f$ induces a chain map $f_{*} : M_{*} \longrightarrow M'_{*}$
and because $f_{*} (z_{V})_{*}$ is chain homotopic to  $(z_{V'})_{*}f_{*}$ the pair of ${\mathcal A}$-morphisms $\nu(f_{i}) \nu(z_{V})_{i}$ and $\nu(z_{V'})_{i} \nu(f_{i})$ for $i=0,1$ induce $f z_{V} = z_{V'}f$ and so $ \nu(z_{V})_{i}$ and $\nu(z_{V'})_{i}$ for $i=0,1$ induce the elements $z_{V}, z_{V'}$ of the central family.

Conversely the degree $0$ and $1$, for any choice of resolution of $V$ determine a central morphism $z_{V}$. When ${\mathcal A}$ is the category of smooth representations of $G$ the morphisms 
$(z_{V})_{i}$ for $i=0,1$ are described in terms of elements of the hyperHecke algebra satisfying the monocentric conditions of Remark \ref{3.6}.
\begin{question}
\label{11.1.4}
\begin{em}

(i) \ Very briefly, the method of \cite{PD84} uses the connection between Hecke modules and smooth representation to reduce determination of the Bernstein centre to the calculation of the centre of the algebra $\lim_{\leftarrow} {\mathcal H}_{G}e$ given by inverse limit over idempotents, $e$, of ${\mathcal H}_{G}$.

In view of the occurrence of the clearly related Hecke and hyperHecke algebras into the topic of this section, I am very curious to further understand the dictionary between the two approaches.

(ii) \  The conditions of Remark \ref{3.6} would be most explicit in the case of $G = GL_{2}K$ where $K$ is a non-Archimedean local field which is given by ``explicit'' homological algebra in (\cite{Sn18} Chapter Two).

Is it possible to describe the dictionary of (i) in this case?
\end{em}
\end{question}
\newpage

\begin{dummy}{A little more details regarding Question \ref{11.1.4}}
\label{11.1.5}
\begin{em}

Possibly the question deserves further elaboration.
\end{em}
\end{dummy}

 Let $G$ be a locally profinite group, such a group is locally compact and totally discontinuous as in (\cite{PD84} p. 2). Recall that the Hecke algebra of a locally compact, totally disconnected group is an  idempotented algebra. Assume that $G$ is unimodular - that is, the left invariant Haar measure equals the right-invariant Haar measure of $G$ (\cite{DB96} p.137).
  
 The Hecke algebra of $G$, denoted by ${\mathcal H}_{G}$ is the space $C_{c}^{\infty}(G)$ of locally constant,compactly supported $k$-valued functions (or measures) on $G$ where $k$ is an algebraically closed field of characteristic zero. In the case of the ``classical Langlands'' situation $k = {\mathbb C}$, the complex numbers.
 
 The algebra structure on ${\mathcal H}_{G}$ is given by the convolution product (\cite{DB96} p.140 and p.255 and \S13)
 \[   (\phi_{1} * \phi_{2})(g)  =  \int_{G} \phi_{1}(gh) \phi_{2}(h^{-1}) dh  = \int_{G}  \phi_{1}(h) \phi_{2}(h^{-1}g) dh  . \]
 This integral requires only one of $\phi_{1}, \phi_{2}$ to be compactly supported in order to land in ${\mathcal H}_{G}$.
 
 Suppose that $K_{0} \subseteq G$ is a compact, open subgroup. Define an idempotent
 \[   e_{K_{0}} = \frac{1}{{\rm vol}(K_{0})} \cdot \chi_{K_{0}} \]
 where $\chi_{K_{0}}$ is the characteristic function of $K_{0}$. If $K_{0} \subseteq K_{1}$ then $e_{K_{0}} * e_{K_{1}} = e_{K_{1}} $. 
 An admissible $k$-representation $V$ of $G$ is equivalent to a non-degenerate ${\mathcal H}_{G}$-module (\cite{PD84} p.2). The dictionary between these concepts may be described in the following manner. Let $v \in V$ be fixed by a compact open subgroup $K$ (i.e. $v \in V^{(K,1)}$). Then we assert that the idempotent $e_{K}$ shall satisfy $e_{K} \cdot v = v$. Then $\pi : G \longrightarrow {\rm Aut}_{k}(V)$ will be defined by $\pi(g)(v) = (\delta_{g} * e_{K})v$ where $(- * -)$ is the convolution product and $\delta_{g})$ is the measure giving rise to evaluation at $g$ on functions (\cite{AJS79} pp.18-22).
 
 ${\mathcal H}_{G}$ is an idempotented algebra (\cite{DB96} p. 309; \cite{PD84} p.2) so that 
  ${\mathcal H}_{G} =  {\mathcal H}_{G}e \oplus  {\mathcal H}_{G}(1-e)$ for any idempotent $e$. The topology of  ${\mathcal H}_{G}$ has a basis consisting of $\{ 0 \}$ and the ${\mathcal H}_{G}(1-e)$'s.
Define $\hat{{\mathcal H}}_{G}$ to be the completion of ${\mathcal H}_{G}$ in this topology. Therefore 
 $\hat{{\mathcal H}}_{G}$ is the inverse limit of the $\hat{{\mathcal H}}_{G}e$'s over the transition maps 
 $\hat{{\mathcal H}}_{G}f  \mapsto \hat{{\mathcal H}}_{G}e$ given by $x \mapsto xe$ where $Hr \leq Hf$ (i.e. $e<< f$ in the poset of idempotents).
 
 Then $\hat{{\mathcal H}}_{G} \subset  {\rm End}_{{\mathbb Z}}({\mathcal H}_{G} )$ (\cite{PD84} p.3) and the centre of $\hat{{\mathcal H}}_{G}$ gives the Bernstein centre of the category of admissible $k$-representations $V$ of $G$ (\cite{PD84} \S\S1.3-1.7).
 
 Now let us consider how to adapt the approach of \ref{11.1.3} to study the Bernstein centre of the category of admissible $k$-representations of $G$ by means of the monomial category $_{k[G], \underline{\phi}}{\bf mon}$. For this purpose I have opted for the category of  admissible $k$-representations of $G$ with a fixed central character $\underline{\phi}$. In the semi-simple situation such as the ``classical Langlands'' this ought not to be a loss of generality. Resolutions of such representations fall into the situation of \ref{11.1.3} in which ${\mathcal B} = _{k[G], \underline{\phi}}{\bf mon}$ and ${\mathcal A} = _{k[G], \underline{\phi}}{\bf mod}$. The construction, in \S\S 8-10, begins with full embeddings of both the monomial category and the module category into a functor category (\ref{8.1.11} and \ref{8.1.13})
 \[ {\mathcal I} :   _{k[G], \underline{\phi}}{\bf mod}  \longrightarrow  funct_{k}^{o}(_{k[G], \underline{\phi}}{\bf mon}, _{k}{\bf mod})   \]
 and
  \[  {\mathcal J} :   _{k[G], \underline{\phi}}{\bf mon}  \longrightarrow  funct_{k}^{o}(_{k[G], \underline{\phi}}{\bf mon}, _{k}{\bf mod}) . \]
  
 The approach of \ref{11.1.3} to describe the Bernstein centre of $ _{k[G], \underline{\phi}}{\bf mod} $
 will be to replace the chain complex $M_{*}$ of \ref{11.1.3} by the monomial resolution of each addmissible $V  \in  _{k[G], \underline{\phi}}{\bf mod} $. Then we attempt to find the $k[G]$-endomorphisms of $V$, which are necessary to give the Bernstein centre of the category, by finding suitable endomorphisms of the monomial resolution. 
 
 This approach will lead us to the $k$-group ring of the monocentre of the hyperhecke algebra, which was described in \S3. As mentioned in Remark \ref{3.6}, it may even suffice to describe a partial chain map endomorphism of the monomial resolution which is defined only on degrees $0$ and $1$ (and possibly $2$ also) since that is enough to define a $k[G]$-endomorphisms of $V$.
 
 The bar-monomial resolution is defined via the full embedding ${\mathcal I} $ by means of which the Bernstein category of the category of admissible representations becomes a Bernstein centre of a subcategory of the functor category $ funct_{k}^{o}(_{k[G], \underline{\phi}}{\bf mon}, _{k}{\bf mod}) $.
 
 As in Theorem \ref{8.1.19} we set $S \in _{k[G], \underline{\phi}}{\bf mon}$ to be 
\[ S = \oplus_{(H, \phi) \in  {\mathcal M}_{\underline{\phi}}(G)} \ \underline{{\rm Ind}}_{H}^{G}(k_{\phi})  \]
which, to be precise, is the direct sum of $ \underline{{\rm Ind}}_{H}^{G}(k_{\phi})$'s - one for each $G$-orbit of $ {\mathcal M}_{\underline{\phi}}(G)$ 
(i.e. one for each $ _{k[G], \underline{\phi}}{\bf mon}$-isomorphism class of irreducible objects). The ring of endomorphisms of $S$ is denoted by ${\mathcal A}_{S}$; it is isomorphic to the hyperHecke algebra.

In Theorem \ref{8.1.19} we have an equivalence  $ \Phi_{S}$ from the functor category to the category of hyperHecke modules. In \S8 the context described assumes that $G$ is finite modulo the centre. However in \S9 it is observed that everything extends to the case where $G$ is compact open modulo the centre.

In addition, the results still work if $S$ is enlarged by replacing the object by its finite direct sum. In this case ${\mathcal A}_{S}$ is changed by a Morita equivalence.

In the situation thus far the bar-monomial resolution of $V$ is transported to the ${\mathcal A}_{S}$-bar resolution of 
\[ {\mathcal I}(V)(S) =  {\rm Hom}_{_{k[G], \underline{\phi}}{\bf mod} }({\mathcal V}(S), V) . \]

The case of a general locally profinite group is treated as follows. $G$ contains a finite number of maximal compact open modulo the centre subgroups, $\{ H_{i} \}$. By choosing suitable $S_{H_{i}}$'s for each of these subgroups we can form a local system, of the associated ${\mathcal A}_{S}$-bar resolutions, on the Tammo to Dieck space $\underline{E}(G, {\mathcal C})$ as in \S10. The resulting bicomplex of ${\mathcal A}_{S}$-modules has a total complex which is the image of the bar-monomial resolution of $V$ in 
$   _{k[G], \underline{\phi}}{\bf mon}   $.

This was the method adopted in \S\S8-10 to construct the bar-monomial resolution of an admissible 
$k[G]$-representation with central character $\underline{\phi}$. However there is clearly a homological chain equivalence between the ${\mathcal A}_{S}$-bar resolution of $ {\mathcal I}(V)(S) $ and the total complex of \S\S8-10. The latter is a chain complex of hyperHecke modules, which  we can use in \ref{11.1.3}. The beginning of the chain complex has the form
\[  0 \leftarrow  {\mathcal I}(V)(S)  \leftarrow  {\mathcal I}(V)(S)  \otimes S \leftarrow 
 {\mathcal I}(V)(S) \otimes {\mathcal A}_{S} \otimes S \leftarrow \ldots .  \]

Clearly an element in the centre of ${\mathcal A}_{S}$, which contains the group-ring over $k$ of the monocentre of $G$, induces an ``Bernstein centre'' endomorphism of this complex which originates from one for $V$.

Finally, there is a modification which we can make by taking the quotient of ${\mathcal A}_{S}$ by killing the 
$ \underline{{\rm Ind}}_{H}^{G}(k_{\phi})$'s for which $V^{(H, \phi)} = 0$. The ``monocentre'' of this quotient might possibly be larger than that of ${\mathcal A}_{S}$ - thereby possibly giving more $k[G]$-endomorphisms of $V$ corresponding to elements in the Bernstein centre of the category of admissible $k[G]$-representations with central character $\underline{\phi}$.

To reiterate Question \ref{11.1.4}:
\begin{question}
\label{11.1.6}
\begin{em}

Question \ref{11.1.4} asks: Does this ``dictionary'' all make sense? If so, bearing in mind the material of \S7 which expresses ${\mathcal A}_{S}$ in terms of convolution products in the classical Hecke algebra, what is the precise dictionary between this section and the $\hat{{\mathcal H}}_{G}$ description of the Bernstein centre given in \cite{PD84}?
\end{em}
\end{question}

 \section{Appendix: Comparison of inductions}
 
 In the case of finite groups this Appendix compares the ``tensor product of modules'' model of an induced representation with the ``function space'' model\footnote{In (\cite{Sn18} Chapter Two, Definition 1.1) my unreliable typography resulted in a superfluous suffix ``-1'' which gives the right action. This this essay I have been more careful to give the correct formula for the left action, since left actions are my usual preference.}. 
 
  Suppose that $H \subseteq G$ are finite groups and that $W$ is a vector space over an algebraically closed field $k$ together with a left $H$-action given by a homomorphism
 \[   \phi : H \longrightarrow  {\rm Aut}_{k}(W)  . \]
 In this case the functional model for the induced representation is given by the $k$-vector space of functions $X_{(H, \phi)}$ consisting of functions of the form
  $f : G \longrightarrow  W $ such that $f(hg) = \phi(h)(f(g))$. The left $G$-action on these functions is given by $(g \cdot f)(x) = f(xg)$ since $((g_{1} g_{2}) \cdot f)(x) = f(xg_{1}g_{2})$ but if $\xi(y) = f(yg_{2})$
then $(g_{1} \cdot \xi)(z) = \xi( z g_{1})  = f( zg_{1} g_{2})$ so that
\[ g_{1} \cdot ( g_{2} \cdot f) =  (g_{1}g_{2}) \cdot f \]
as required.
  
For $w \in W $ we have a function $f_{w}$, supported in $H$ and satisfying 
\linebreak
$(h \cdot f_{w}) = f_{\phi(h)(w)}$ for $h \in H$ so that $f_{w}(1) = w$. We have a left $k[H]$-module map 
\[   f: W \longrightarrow X_{(H, \phi)  } \]
defined by $w \mapsto f_{w}$. This makes sense because
\[ \begin{array}{l}
f_{w}(x) = \left\{  \begin{array}{ll}
\phi(x)(w) & {\rm if} \ x \in H , \\
\\
0 & {\rm if} \ x \not\in H
\end{array} \right.
\end{array} \]
so that
\[ \begin{array}{l}
f_{\phi(h)(w)}(x) = \left\{  \begin{array}{ll}
\phi(x)( \phi(h)(w) )  = \phi(xh)(w)& {\rm if} \ x \in H , \\
\\
0 & {\rm if} \ x \not\in H
\end{array} \right.   .
\end{array} \]
The map $f$ induces a left $k[G]$-module map, which is an isomorphism,
\[       \hat{f} : Ind_{H}^{G}( W ) = k[G]  \otimes_{k[H]} W  \stackrel{\cong}{\longrightarrow}   X_{(H, \phi)  }    \]
given by $\hat{f}(g \otimes_{k[H]} w) = g \cdot f_{w}$. This is well-defined, as is well-known, because
\[  \hat{f}(g \otimes \phi(h)(w))  = g \cdot f_{\phi(h)(w)} = g \cdot ( h \cdot f_{w} ) = 
\hat{f}(gh \otimes_{k[H]} w) .  \]

Henceforth, in this Appendix, I shall consider only the case when 
\linebreak
${\rm dim}_{k}(W) =1$. In this case $W = k_{\phi}$ will denote the $H$-representation given by the action  
$h \cdot v = \phi(h)v$ for $h \in H, v \in k$.

 As in Definition \S2,  write 
$[(K, \psi), g, (H, \phi)]$ for any triple consisting of $g \in G$, characters $\phi, \psi$ on 
subgroups $H, K \leq G$, respectively such that 
\[    (K, \psi) \leq (g^{-1}Hg, (g)^{*}(\phi)) \]
which means that $K \leq  g^{-1}Hg$ and that $\psi(k) = \phi(h)$ where $k = g^{-1}hg$ for $h \in H, k \in K$.

We have a well-defined left  $k[G]$-module homomorphism
\[    [(K, \psi), g, (H, \phi)] : k[G] \otimes_{k[K]} k_{\psi}  \longrightarrow  k[G] \otimes_{k[H]} k_{\phi}  \]
given by the formula $[(K, \psi), g, (H, \phi)](g' \otimes_{k[K]} v) = g' g^{-1} \otimes_{k[H]} v$.

This is well-defined because, for $k = g^{-1}hg \in K$ we have
\[   \begin{array}{ll}
[(K, \psi), g, (H, \phi)](g'  k \otimes_{k[K]} \psi(k)^{-1}(v)) &  = 
g'  k g^{-1} \otimes_{k[H]} \psi(k)^{-1}(v)   \\
\\
& =  g'  g^{-1} (gk g^{-1}) \otimes_{k[H]} \psi(k)^{-1}(v) \\
\\
& =   g'  g^{-1}  \otimes_{k[H]}  \phi(h) (\psi(k)^{-1}(v))  \\
\\
&  =  g'  g^{-1}  \otimes_{k[H]}  v    \\
\\
& = [(K, \psi), g, (H, \phi)](g'   \otimes_{k[K]} v) . 
\end{array} \]
Being well-defined, it is clearly a  left  $k[G]$-module homomorphism. 

When two of these triples correspond to $k[G]$-modules homomorphisms which are composeable we have
\[    \begin{array}{l}
[(H, \phi), g_{1}, (J, \mu)][(K, \psi), g, (H, \phi)] (g' \otimes_{k[K]} v)   \\
\\
=  [(H, \phi), g_{1}, (J, \mu)] (g'  g^{-1} \otimes_{k[H]} v) \\
\\
=  g'  g^{-1} g_{1}^{-1} \otimes_{k[J ]} v)   \\
\\
=  g'  (g_{1}g)^{-1} \otimes_{k[J ]} v) \\
\\
= [(K, \phi), g_{1}g, (J, \mu)] (g' \otimes_{k[K]} v) ,
\end{array}   \]
which makes sense because 
\[   (K, \psi) \leq (g^{-1}Hg, (g)^{*}(\phi))  \leq   (g^{-1} g_{1}^{-1}J g_{1}g, (g_{1}g)^{*}(\phi)) .\]

In order to define a left $k[G]$-homomorphism
\[     [(K, \psi), g, (H, \phi)] :    X_{(K, \psi)} \longrightarrow  X_{(H, \phi)} \]
satisfying the relation 
\[ \hat{f} \cdot  [(K, \psi), g, (H, \phi)] =  [(K, \psi), g, (H, \phi)]  \cdot \hat{f} :
   k[G] \otimes_{k[K]} k_{\psi}   \longrightarrow  X_{(H, \phi)} \]
  we set
  \[    [(K, \psi), g, (H, \phi)] (  g_{1} \cdot f_{v}) = (g_{1}g^{-1}) \cdot f_{v}  . \]
  
  This is well-defined because $g_{1} \cdot f_{v} = g_{2} \cdot f_{v} \in X_{(K, \psi)}$ if and only if  
  \linebreak
  $(g_{1} \cdot f_{v})(x) = (g_{2} \cdot f_{v})(x)$ for all $x \in G$ and so, for $v \not= 0$,
\[  f_{v}(xg_{1}) =  \left\{  \begin{array}{ll}
\psi(xg_{1})(v) & {\rm if} \ xg_{1}  \in K \\
\\
0 &  {\rm if} \ xg_{1} \not\in K \\
\end{array} \right.  \]
equals
\[  f_{v}(xg_{2}) =  \left\{  \begin{array}{ll}
\psi(xg_{2})(v) & {\rm if} \ xg_{2} \in K \\
\\
0 &  {\rm if} \ xg_{2}^{-1} \not\in K  . \\
\end{array} \right.  \]
Hence $k =  (xg_{1})^{-1} xg_{2} = g_{1}^{-1} g_{2} \in K$ and
\[  \psi(xg_{1})(v) =   \psi(xg_{2})(v) = \psi(xg_{1})(  \psi(g_{1}^{-1}g_{2})(v)) \]
so that $ \psi(g_{1}^{-1}g_{2})(v)  =  v$. 

Conjugating by $g$ we have $gkg^{-1} = g g_{1}^{-1} g_{2} g^{-1} \in H$ which implies that
\[  (x g_{1} g^{-1})^{-1} xg_{2}g^{-1} = g g_{1}^{-1}  x^{-1} x g_{2}g^{-1} = 
g g_{1}^{-1} g_{2}g^{-1}  \in H  \]
and so $xg_{1}g^{-1} \in H$ if and only if
$ x g_{1}g^{-1} g g_{1}^{-1} g_{2}^{-1} g^{-1} =  x  g_{2}^{-1} g^{-1}  \in H $. Also 
\[     \begin{array}{ll}
 \phi(xg_{2}g^{-1})(v)  & =  \phi( xg_{1}g^{-1} g g_{1}^{-1} g_{2} g^{-1})(v) \\
 \\
 &  = \phi(xg_{1}g^{-1})( \phi(g g_{1}^{-1} g_{2} g^{-1})(v)  )  \\
 \\
 & =   \phi(xg_{1}g^{-1})( \psi(g_{1}^{-1} g_{2} )(v)  )  \\
 \\
 & =  \phi(xg_{1}g^{-1})(v)  .
 \end{array}  \]
 Therefore
 \[       [(K, \psi), g, (H, \phi)] (  g_{1} \cdot f_{v}) =   [(K, \psi), g, (H, \phi)] (  g_{2} \cdot f_{v}) ,  \]
 as required in order for the left $k[G]$-module homomorphism to be well-defined.
 
 It is easy to see that transporting the map $[(K, \psi), g, (H, \phi)]$ from the tensor product model of the induced representation to the function space model gives the left $k[G]$-homomorphism 
 whose well-definedness we have just verified.

Among the left $k[G]$-maps
\[  k[G] \otimes_{k[K]} k_{\psi}  \longrightarrow  k[G] \otimes_{k[H]} k_{\phi}  \]
we have the relations, $h \in H, k \in K$
 \[   [(K, \psi), gk, (H, \phi)]  =  [(K, \psi), g, (H, \phi)] \cdot (1 \otimes_{k[K]} \psi(k^{-1})) \]
and
 \[     [(K, \psi), hg, (H, \phi)]  =  (1 \otimes_{k[H]} \phi(h^{-1})) \cdot  [(K, \psi), g, (H, \phi)]   . \]
This is  because 
\[  \begin{array}{ll}
g_{1} (gk)^{-1} \otimes_{k[H]} v   & =  g_{1} g^{-1} gk^{-1}g^{-1} \otimes_{k[H]} v  \\
\\
& = g_{1} g^{-1}  \otimes_{k[H]} \phi(gk^{-1}g^{-1})(v) \\
\\
& = g_{1} g^{-1}  \otimes_{k[H]} \psi(k^{-1})(v) 
\end{array}   \]
and 
\[  g_{1} (hg)^{-1}  \otimes_{k[H]} v = g_{1} g^{-1} h^{-1} \otimes_{k[H]} v
 g_{1} g^{-1}  \otimes_{k[H]} \phi(h^{-1})v.   \]

Transporting these relations by $\hat{f}$ yields
 \[   [(K, \psi), gk, (H, \phi)]  =  [(K, \psi), g, (H, \phi)]  \cdot (g \cdot f_{v} \mapsto  g \cdot  f_{\psi(k^{-1})(v)})  \]
and
 \[     [(K, \psi), hg, (H, \phi)]  = (g \cdot f_{v} \mapsto  g \cdot  f_{\psi(h^{-1})(v)}) \cdot [(K, \psi), g, (H, \phi)]     \]
 because, for example, $g \otimes_{k[K]} v \longrightarrow  g \otimes_{k[K]} \psi(k^{-1})(v)$ is transported to 
 \[  g \cdot f_{v} \mapsto  g \cdot  f_{\psi(k^{-1})(v)} = gk^{-1} \cdot f_{v}  . \]

 \begin{dummy}
 \label{3.1}
 \begin{em}
 
 Now we shall transport the Double Coset Formula over to the functional models.
 \end{em}
 \end{dummy}

The Double Coset Formula (\cite{Sn94} Theorem 1.2.40)is an isomorphism describing the restriction of an induced representation. It is a consequence of the $J$-orbit structure of the left action of a subgroup $J \subseteq G$ on $G/H$.
This is a left $k[J]$-isomorphism of the form 
\[  {\rm Res}_{J}^{G} {\rm Ind}_{H}^{G}(k_{\phi})  \stackrel{\alpha}{\longrightarrow}
\oplus_{z \in J \backslash G  / H} \ {\rm Ind}_{J \bigcap zHz^{-1}}^{J}((z^{-1})^{*}(k_{\phi}))  \]
given by $\alpha(g \otimes_{H} v) = j \otimes_{J \bigcap zHz^{-1}} \phi(h)(v)$ for $g = jzh, j \in J, h \in H$.
The inverse of $\alpha$ is given by $\alpha^{-1}(j \otimes_{J \bigcap zHz^{-1}} v) =
 jz \otimes_{H} (v)$. This makes sense for $\alpha$ because if, for example, $jzh = j' z h'$ then $(j')^{-1} j = z h' h^{-1} z^{-1} \in J \bigcap zHz^{-1}$ so that
\[   \begin{array}{l}
 j' \otimes_{J \bigcap zHz^{-1}} \phi(h')(v)  \\
 \\
 =    j'  (j')^{-1} j \otimes_{J \bigcap zHz^{-1}} \phi( h' h^{-1}  )^{-1}(\phi(h')(v)) \\
 \\
 =    j \otimes_{J \bigcap zHz^{-1}}  \phi(h)(v)
 \end{array}  \]
 and for $\alpha^{-1}$,
  \[ \begin{array}{ll}
 \alpha^{-1}( \alpha(  g \otimes_{H} v))   &
 = \alpha^{-1}( j \otimes_{J \bigcap zHz^{-1}} \phi(h)(v)) \\
 \\
 &  =   jz \otimes_{k[H]} \phi(h)(v) \\
  \\
 &  =  g \otimes_{H} v  .
 \end{array} \]

 In terms of functional models for induced representations we have a left $k[J]$-module isomorphism,
obtained by transporting the above isomorphism via the appropriate map $\hat{f}$ introduced above, of the form
 \[   a:  {\rm Res}_{J}^{G} X_{(H, \phi)}  
 \stackrel{\cong}{\longrightarrow}  
 \oplus_{z \in J \backslash G  / H} \   X_{ (J \bigcap zHz^{-1}, (z^{-1})^{*}(\phi)) } \]
 whose $ X_{ (J \bigcap zHz^{-1}, (z^{-1})^{*}(\phi)) }$-component equals
 \[    a(g \cdot f_{v}) =  j \cdot f_{\phi(h)( v) }  \in X_{ (J \bigcap zHz^{-1}, (z^{-1})^{*}(\phi)) } \ {\rm if } \ g = jzh . \]
 The inverse, $a^{-1}$, maps $j \cdot f_{v} \in X_{ (J \bigcap zHz^{-1}, (z^{-1})^{*}(\phi)) } $ by the formula
 \[   a^{-1}(j  \cdot f_{v})  =  jz  \cdot f_{v}  \in X_{(H, \phi)} . \]
 Thus $a$ is clearly a left $k[J]$-module map.

\section{Appendix: Smooth representations and Hecke modules}

In this Appendix, for my convenience, representations are complex representations.

Now let $\Gamma$ be a compact totally disconnected group. Denote by $\hat{\Gamma}$ the set of equivalence classes of finite-dimensional irreducible representations of $\Gamma$ whose kernel is open - and hence of finite index in $\Gamma$.

Suppose now that $\Gamma$ is finite  and $(\pi, V)$ is a representation of $\Gamma$ on a possible infinite dimensional vector space $V$. If $\rho \in \hat{\Gamma}$ let $V(\rho)$ be the sum of all invariant subspaces of $V$ that are isomorphic as $\Gamma$-modules to $V_{\rho}$. $V(\rho)$ is the $\rho$-isotypic subspace of $V$. We have
\[ V \cong \oplus_{\rho \in \hat{\Gamma}} \ V_{\rho} . \]

Now we generalise this to smooth representations of a totally disconnected locally compact group. Choose a compact open subgroup $K$ of $G$. The compact open normal subgroups of $K$ form a basis of neighbourhoods of the identity in $K$. Let $\rho  \in \hat{K}$ then the kernel of $\rho$ is $K_{\rho}$ a compact open normal subgroup of finite index.
\begin{proposition}{(\cite{DB96} Proposition 4.2.2)}
\label{12.2}
\begin{em}

Let $(\pi, V)$ be a smooth representation of $G$.
Then
\[ V \cong \oplus_{\rho \in \hat{K}} \ V_{\rho} . \]
The representation $\pi$ is admissible if and only if each $V(\rho)$ is finite-dimensional.
\end{em}
\end{proposition}

Let $(\pi, V)$ be a smooth representation of $G$. If $\hat{v} : V \longrightarrow  {\mathbb C}$ is a linear functional we write $\langle v , \hat{v} \rangle = \hat{v}(v)$ for $v \in V$. We say $\hat{v}$ is smooth if there exists an open neighbourhood $U$ of $1 \in G$ such that for all $g \in U$
\[  \langle \pi(g)(v) , \hat{v} \rangle = \hat{v}(v) . \]
Let $\hat{V}$ denote the space of smooth linear functionals on $V$.

Define the contragredient representation $(\hat{\pi}, \hat{V})$ is defined by
\[   \langle v ,   \hat{\pi}(g)(\hat{v}) \rangle =   \langle \pi(g^{-1})(v) , \hat{v} \rangle . \]
The contragredient representation of a smooth representation is a smooth representation. Also
\[ \hat{V} \cong \oplus_{\rho \in \hat{K}} \ V_{\rho}^{*}  \]
where $V_{\rho}^{*}$ is the dual space of $V_{\rho}$.

Since the dual of a finite-dimensional $V_{\rho}$ is again finite-dimensional the contragredient of an admissible representation is also admissible. Also $\hat{\hat{\pi}} = \pi$.

If $X$ is a totally disconnected space a complex valued function $f$ on $X$ is smooth if it is locally constant.
Let ${\mathcal H}_{G}$ be, as before,  the space of smooth compactly supported complex-valued functions on $X=G$. Assuming $G$ is unimodular ${\mathcal H}_{G}$ is an algebra without unit under the convolution product
\[    (\phi_{1} * \phi_{2})(g) = \int_{G} \ \phi_{1}(gh^{-1}) \phi_{2})(h) dh .  \]
This is the Hecke algebra - an idempotented algebra (see \S6).

If $\phi \in {\mathcal H}$ define $\pi(\phi) \in {\rm End}(V)$ with $V$ as above 
\[   \pi(\phi)(v) = \int_{G} \phi(g) \pi(g)(v) dg  . \]
Then 
\[    \pi (\phi_{1} * \phi_{2})  =  \pi(\phi_{1}) \cdot \pi(\phi_{2})  \]
so that $V$ is an ${\mathcal H}$-representation.

The integral defining $\phi$ may be replaced by a finite sum as follows. Choose an open subgroup $K_{0}$ fixing $v$. Choosing $K_{0}$ small enough we may assume that the support of $\phi$ is contained in a finite union of left cosets $\{ g_{i}K_{0} \ | \ 1 \leq i \leq t \}$. Then
\[   \pi(\phi)(v) =  \frac{1}{{\rm vol}(K_{0})}  \sum_{i=1}^{t}  \phi(g_{i})  \pi(g_{i})(v) . \]

{\bf Finite group example:}

Let $(\pi, V)$ be a finite-dimensional representation of a finite group $G$. Write ${\mathcal H}$ for the space of functions from $G$ to ${\mathbb C}$. If $\phi_{1}, \phi_{2} \in {\mathcal H}$ define $\phi_{1} * \phi_{2} 
\in {\mathcal H}$ by
\[   (\phi_{1} * \phi_{2})(g) = \sum_{h \in G} \ \phi_{1}(gh^{-1}) \phi_{2}(h) . \]
For $\phi \in {\mathcal H}$ define $\pi(\phi) \in {\rm End}_{{\mathbb C}}(V)$ by
\[   \pi(\phi)(v) = \sum_{g \in G} \ \phi(g) \pi(g)(v) . \]
Hence
\[ \begin{array}{l}
\pi(\phi_{1}(  \pi(\phi_{2})(v) )   \\
\\
=  \pi(\phi_{1})(   \sum_{g \in G} \ \phi_{2}(g) \pi(g)(v)  ) \\
\\
=   \sum_{g \in G} \ \phi_{2}(g)  \pi(\phi_{1}( \pi(g)(v))  \\
\\
=    \sum_{g \in G} \   \phi_{2}(g)  \sum_{\tilde{g} \in G} \ \phi_{1}(\tilde{g}) (\pi(\tilde{g}( \pi(g)(v)) \\
\\
=   \sum_{g, \tilde{g} \in G} \  \phi_{2}(g)  \phi_{1}(\tilde{g}) (\pi(\tilde{g}g)(v)) .
\end{array} \]

Now 
\[ \begin{array}{l}
\pi( \phi_{1} * \phi_{2})(v) \\
\\
= \sum_{g_{1} \in G} \ ( \phi_{1} * \phi_{2})(g_{1}) \pi(g_{1})(v)  \\
\\
= \sum_{g_{1}, h \in G} \  \phi_{1}(h_{1}h^{-1}) \phi_{2}(h) \pi(g_{1})(v)  .
\end{array} \]
Setting $g=h$, $\tilde{g}g = g_{1}$ shows that 
\[  \pi( \phi_{1} * \phi_{2}) =  \pi(\phi_{1} \cdot \pi(\phi_{2}) .\]

Also ${\mathcal H} \cong {\mathbb C}[G]$  because if $f_{g}(x) = 0$ if $g \not= x$ and $f_{g}(g)=1$ then 
\[ f_{g} * f_{g'} = f_{gg'} . \]

\begin{proposition}{(\cite{DB96} Proposition 4.2.3)}
\label{12.3}
\begin{em}

Let $(\pi, V)$ be a smooth non-zero representation of $G$. Then equivalent are:

(i)  \ $\pi$ is irreducible.

(ii) \ $V$ is a simple ${\mathcal H}$-module.

(iii) \ $V^{K_{0}}$ is either zero or simple as an ${\mathcal H}_{K_{0}}$-module for all open subgroups $K_{0}$.
Here ${\mathcal H}_{K_{0}} = e_{K_{0}} * {\mathcal H}  * e_{K_{0}}$.
\end{em}
\end{proposition}

Schur's Lemma holds (\cite{DB96} \S4.2.4)for $(\pi, V)$ an irreducible admissible representation of a totally disconnected group $G$.
\begin{proposition}{(\cite{DB96} Proposition 4.2.5)}
\label{12.4}
\begin{em}

Let $(\pi, V)$ be an admissible representation of the totally disconnected locally compact group $G$ with contragredient $(\hat{\pi}, \hat{V})$. Let $K_{0} \subseteq G$ be a compact open subgroup. Then the canonical pairing between $V$ and $\hat{V}$ induces a non-degenerate pairing betweem $V^{K_{0}}$ and 
$\hat{V}^{K_{0}}$.
\end{em}
\end{proposition}

{\bf The trace}

As with representations of finite groups the character of an admissible representation of a totally disconnected locally compact group $G$ is an important invariant. It is a distribution. It is a theorem of Harish-Chandra that if $G$ is a reductive $p$-adic group then the character is in fact a locally integrable function defined on a dense subset of $G$.

We shall define the character as a distribution on ${\mathcal H}_{G} = C_{c}^{\infty}(G)$. Suppose that $U$ is a finite-dimensional vector space and let $f:U \longrightarrow U$ be a linear map. Suppose ${\rm Im}(f) \subseteq U_{0} \subseteq U$. Then we have 
\[ {\rm Trace}(f:U_{0} \longrightarrow U_{0}) = {\rm Trace}(f:U \longrightarrow U) . \]
Therefore we may define the trace of any endomorphism $f$ of $V$ which has finite rank by choosing any
finite-dimensional $U_{0}$ such that ${\rm Im}(f) \subseteq U_{0} \subseteq V$ and by defining
\[   {\rm Trace}(f) =   {\rm Trace}(f:U_{0} \longrightarrow U_{0}) . \]

Now let $(\pi, V)$ be an admissible representation of $G$. Let $\phi \in {\mathcal H}_{G}$. Since $\phi$ is compactly supported and locally constant there exists a compact open $K_{0}$ such that $\phi \in {\mathcal H}_{K_{0}}$. The endomorphism $\pi(\phi)$ has image in $V^{K_{0}}$ which is finite-dimensional - by admissibility - so we define the trace distribution
\[ \chi_{V} : {\mathcal H} \longrightarrow {\mathbb C} \]
by
\[   \chi_{V}(\phi) = {\rm Trace}(\pi(\phi)) . \]
\begin{proposition}{(\cite{DB96} Proposition 4.2.6)}
\label{12.5}
\begin{em}

Let $R$ be an algebra over a field $k$. Let $E_{1}$ and $E_{2}$ be simple $R$-modules that are finite-dimensional over $k$. For each $\phi \in R$ if 
\[ {\rm Trace}( (\phi \cdot -):E_{1} \longrightarrow E_{1}) =  {\rm Trace}( (\phi \cdot -):E_{2} \longrightarrow E_{2})   \]
then the $E_{i}$ are isomorphic $R$-modules.
\end{em}
\end{proposition} 
\begin{proposition}{(\cite{DB96} Proposition 4.2.7)}
\label{12.6}
\begin{em}

Let $(\pi_{1}, V_{1})$ and $(\pi_{2}, V_{2})$ be irreducible admissible representations of $G$ (as above) such that, for each compact open $K_{1}$, $V_{1}^{K_{1}} \cong V_{2}^{K_{1}}$ as ${\mathcal H}_{K_{1}}$-modules then $(\pi_{1}, V_{1}) \cong (\pi_{2}, V_{2})$.
\end{em}
\end{proposition} 
\begin{theorem}{(\cite{DB96} Theorem 4.2.1)}
\label{12.7}
\begin{em}

Let $(\pi_{1}, V_{1})$ and $(\pi_{2}, V_{2})$ be irreducible admissible representations of $G$ (as above) such that $\chi_{V_{1}} = \chi_{V_{2}}$ then $(\pi_{1}, V_{1}) \cong (\pi_{2}, V_{2})$.
\end{em}
\end{theorem} 

From this one sees that the contragredient of an admissible irreducible $(\pi, V)$ of $GL_{n}K$ ($K$ a $p$-adic local field) 
is given by $\pi_{1}(g) = \pi( (g^{-1})^{tr}) $
on the same vector space $V$.

\end{document}